\documentclass[a4paper,12pt]{article}
\usepackage{epsfig,amssymb,wrapfig,makeidx, fullpage}

\newtheorem{thm}{Theorem}[section]

\newtheorem{cor}[thm]{Corollary}

\newtheorem{lem}[thm]{Lemma}
\newtheorem{propn}[thm]{Proposition}

\newenvironment{ass}[1][Assumption 1]{\begin{trivlist}\item[\hskip\labelsep{\bfseries #1}]}{\end{trivlist}}
\newenvironment{proof}{\\{\bf Proof:}}{\hspace*{\fill}$\square$\endlist\par\bigskip}
\newenvironment{proof11}{\\{\bf Proof of Theorem 1.1:}}{\hspace*{\fill}$\square$\endlist\par\bigskip}
\newenvironment{proof12}{\\{\bf Proof of Theorem 1.2:}}{\hspace*{\fill}$\square$\endlist\par\bigskip}

\begin{document}
\title{Convergence of simple random walks on random discrete trees to Brownian motion on the continuum random tree}
\author{David Croydon\footnote{Dept of Statistics, University of Warwick, Coventry, CV4 7AL, UK; \underline{d.a.croydon@warwick.ac.uk.}}\\ \tiny{UNIVERSITY OF WARWICK}}
\maketitle

\begin{abstract}
In this article it is shown that the Brownian motion on the continuum random tree is the scaling limit of the simple random walks on any family of discrete $n$-vertex ordered graph trees whose search-depth functions converge to the Brownian excursion as $n\rightarrow\infty$. We prove both a quenched version (for typical realisations of the trees) and an annealed version (averaged over all realisations of the trees) of our main result. The assumptions of the article cover the important example of simple random walks on the trees generated by the Galton-Watson branching process, conditioned on the total population size.

\end{abstract}

\section{Introduction}

The goal of this investigation is to provide a description for the scaling limit of the simple random walks on a wide collection of random graph trees. In particular, we will be interested in ordered graph trees whose scaling limit is the continuum random tree of Aldous, see \cite{Aldous2}, and we shall demonstrate that the scaling limit of the associated simple random walks is the Brownian motion on the continuum random tree. This limiting diffusion process was first constructed on typical realisations of the continuum random tree in \cite{Krebs}.

This work falls into the area of random walks in random environments, of which one of the most difficult and interesting examples is the random walk on a critical percolation cluster. One motivation for studying this process is to gain further insight into the conductivity properties of the cluster, about which few rigorous results are known. There is growing evidence that the incipient infinite cluster in high dimensions behaves like the integrated super-Brownian excursion, which may be viewed as the continuum random tree embedded into Euclidean space, see \cite{HaraSlade} and \cite{JM}. Hence investigating the current problem and other properties of the Brownian motion on the continuum random tree may help to increase understanding of this more challenging model.

More immediate implications are provided by the relationship between the continuum random tree and various random graph trees, important examples of which are Galton-Watson process family trees, started from a single ancestor, conditioned on the total population size being $n$. Under the assumptions of a critical, finite variance, non-lattice offspring distribution, it is known that these trees converge to the continuum random tree, see \cite{Aldous3}. This family of trees, with offspring distribution chosen suitably, also provides a representation of a range of combinatorial random trees; a more detailed discussion of such connections is presented in \cite{Aldous2}. The results here provide a rigorous description of the asymptotics of the simple random walks on these sets. A related model is the branching process conditioned to never become extinct, and the transition densities of the simple random walks on these trees were estimated in \cite{BarKum}. It is known that these sets converge to the self-similar continuum random tree, \cite{Aldous2}, and techniques similar to those applied in this article should yield analogous convergence results for these processes.

Both random walks on the incipient infinite percolation cluster and on a critical branching process conditioned to never become extinct were considered by Kesten in \cite{Kesten}. The second of these problems is particularly closely related to ours, and Kesten demonstrates in \cite{Kestenunpub} that the height (distance from the initial ancestor) of the simple random walk on the branching process studied there converges, when rescaled, to a non-trivial limit. Unfortunately, the argument there is long, as complicated branching process arguments were necessary to complete the proof. In essence, the structure of the argument here does owe a debt to this work of Kesten, but by using the ideas provided by Aldous in \cite{Aldous3} for representing abstract trees, we are able to greatly improve the techniques involved and, in the process, generalise the argument, entirely eliminating the need for any branching process arguments. One further advantage we have is knowledge of the limiting set and process, Brownian motion on the continuum random tree. As noted above, the almost-sure existence of this process was initially demonstrated in \cite{Krebs}, but a more concise construction is given in \cite{Croydoncrt}. Using basic properties of this process, and looking at its restriction to finite length sub-trees of the continuum random tree, we are able to employ a ``meet in the middle'' approach for demonstrating our main convergence result, which proves the conjecture of Aldous in \cite{Aldous2}, Section 5.1. We expect that the problem of extending the results proved here to showing that the simple random walk on a critical branching process conditioned to never become extinct converges when rescaled to a related limiting diffusion is merely technical, and may be solved by applying the ideas used here to an increasing sequence of compact subsets of the infinite tree.

To prove our main results, we will work within the framework developed by Aldous in \cite{Aldous1} for building trees as subsets of the Banach space of infinite sequences of real numbers, $l^1$. Throughout, the usual norm on $l^1$ will be denoted by $\|\cdot\|$. We will frequently consider triples of the form $(K,\nu,\mathbf{Q})$, where $K$ is a compact metric space (or finite graph), $\nu$ is a Borel probability measure on $K$ (or a probability measure on the vertices of $K$), and $\mathbf{Q}$ is a probability measure on $C([0,R],K)$ for some $R>0$ (or a probability law on the space of \{0,1,\dots,R\}-indexed processes taking values in the vertices of $K$, respectively). We will say that $(\tilde{K},\tilde{\nu},\tilde{\mathbf{{Q}}})$ is an {\it (isometric) embedding} of $(K,\nu,\mathbf{Q})$ into $l^1$ if there exists a distance-preserving map $\psi:K\rightarrow l^1$ such that $\tilde{K}=\psi(K)$, $\tilde{\nu}=\nu\circ\psi^{-1}$ and $\tilde{\mathbf{Q}}=\mathbf{Q}\circ\psi^{-1}$. In the discrete case, we extend $\tilde{\mathbf{Q}}$ to a probability law on $C([0,R],l^1)$ by linear interpolation of discrete time processes. Note that the triple $(\tilde{K},\tilde{\nu},\tilde{\mathbf{{Q}}})$ is an element of $\mathcal{K}(l^1)\times \mathcal{M}_1(l^1)\times\mathcal{M}_1(C([0,R],l^{1}))$, where $\mathcal{K}(l^1)$ is the space of compact subsets of $l^1$, $\mathcal{M}_1(l^1)$ is the space of Borel probability measures on $l^1$, and $\mathcal{M}_1(C([0,R],l^{1}))$ is the space of Borel probability measures on $C([0,R],l^{1})$. In statements of convergence and distributional results, we assume that the first of these spaces is endowed with the usual Hausdorff topology for compact subsets of $l^1$, and the remaining two are endowed with the topologies induced by the relevant weak convergence. The rescaling operators we will apply to elements of the form $(\tilde{K},\tilde{\nu},\tilde{\mathbf{{Q}}})\in\mathcal{K}(l^1)\times \mathcal{M}_1(l^1)\times\mathcal{M}_1(C([0,R],l^{1}))$ with $R\geq n^{3/2}$ are defined by
\[\Theta_n\left(\tilde{K},\tilde{\nu},\tilde{\mathbf{Q}}\right):=\left(n^{-1/2}\tilde{K},\tilde{\nu}(n^{1/2}\cdot),\tilde{\mathbf{Q}}(\{f\in C([0,R],l^1):\:(n^{-1/2}f(tn^{3/2}))_{t\in[0,1]}\in\cdot\})\right),\]
the images of which are contained in $\mathcal{K}(l^1)\times \mathcal{M}_1(l^1)\times\mathcal{M}_1(C([0,1],l^{1}))$.

The main result of this article is the quenched limit that we prove as Theorem \ref{mainresult}. It describes how, if we have a collection of (deterministic) ordered graph trees $(\mathcal{T}_n)_{n\geq 1}$, where we always assume that $\mathcal{T}_n$ has $n$-vertices, whose search depth functions, $(w_n)_{n\geq 1}$ say, converge when rescaled to a typical realisation of the normalised Brownian excursion, $w$ say, then we can describe precisely the scaling limit of the triple $(\tilde{\mathcal{T}_n},\tilde{\mu}_n,\tilde{\mathbf{P}}^{\mathcal{T}_n}_\rho)$, which is a specific isometric embedding of $({\mathcal{T}_n},{\mu}_n,\mathbf{P}^{\mathcal{T}_n}_\rho)$ into $l^1$, where $\mu_n$ is the uniform measure on the vertices of $\mathcal{T}_n$, and $\mathbf{P}^{\mathcal{T}_n}_\rho$ is the law of the discrete time simple random walk on $\mathcal{T}_n$, started from the root, $\rho=\rho(\mathcal{T}_n)$, of $\mathcal{T}_n$. It is the family of operators $(\Theta_n)_{n\geq 1}$ that we apply to obtain a non-trivial scaling limit, $(\tilde{\mathcal{T}},\tilde{\mu},\tilde{\mathbf{P}}^\mathcal{T}_\rho)$, which is a specific isometric embedding of the triple $(\mathcal{T}_w,{\mu}_w,\mathbf{P}^{\mathcal{T}_w,\mu_w}_\rho)$ into $l^1$. Here, $\mathcal{T}_w$ is the rooted real tree associated with the excursion $w$, (see Section \ref{abstracttree} for an exact definition), $\mu_w$ is the natural measure on $\mathcal{T}_w$, (see (\ref{muw})), and $\mathbf{P}^{\mathcal{T}_w,\mu_w}_\rho$ is the law of the Brownian motion on $(\mathcal{T}_w,\mu_w)$ started from the root $\rho=\rho(\mathcal{T}_w)$, (see Section \ref{abstractproc}).

In the statement of the following result, we assume that $W$ is the normalised Brownian excursion, built on an underlying probability space with probability measure $\mathbf{P}$. Furthermore, we introduce a set $\mathcal{W}^*\subseteq C([0,1],\mathbb{R}_+)$ that satisfies $\mathbf{P}(W\in\mathcal{W}^*)=1$, and which may therefore be thought of as a collection of typical realisations of $W$. The precise definition of $\mathcal{W}^*$ is given at (\ref{wstar}), and a detailed description of the properties of $(\mathcal{T}_w,{\mu}_w,\mathbf{P}^{\mathcal{T}_w,\mu_w}_\rho)$ that hold for  $w\in\mathcal{W}^*$ is given by Lemma \ref{crtprops}.

{\thm \label{mainresult} There exists a set $\mathcal{W}^*\subseteq C([0,1],\mathbb{R}_+)$ with $\mathbf{P}(W\in\mathcal{W}^*)=1$ such that if $(\mathcal{T}_n)_{n\geq 1}$ is a sequence of ordered graph trees whose search-depth functions $(w_n)_{n\geq 1}$ satisfy
\[n^{-1/2}w_n\rightarrow w\]
in $C([0,1],\mathbb{R}_+)$ for some $w\in\mathcal{W}^*$, then there exists, for each $n$, an isometric embedding $(\tilde{\mathcal{T}}_n,\tilde{\mu}_n,\tilde{\mathbf{P}}^{\mathcal{T}_n}_\rho)$ of the triple $({\mathcal{T}}_n,{\mu}_n,{\mathbf{P}}^{\mathcal{T}_n}_\rho)$ into $l^1$ such that
\[\Theta_n\left(\tilde{\mathcal{T}}_n,\tilde{\mu}_n,\tilde{\mathbf{P}}^{\mathcal{T}_n}_\rho\right)\rightarrow(\tilde{\mathcal{T}},\tilde{\mu},\tilde{\mathbf{P}}^\mathcal{T}_\rho)\]
in the space $\mathcal{K}(l^1)\times\mathcal{M}_1(l^1)\times \mathcal{M}_1(C([0,1],l^1))$, where $(\tilde{\mathcal{T}},\tilde{\mu},\tilde{\mathbf{P}}^\mathcal{T}_\rho)$ is an isometric embedding of the triple $(\mathcal{T}_w,{\mu}_w,{\mathbf{P}}^{\mathcal{T}_w,\mu_w}_\rho)$ into $l^1$.}
\bigskip

The choice of embedding of $(\mathcal{T}_w,{\mu}_w,{\mathbf{P}}^{\mathcal{T}_w,\mu_w}_\rho)$ that we use in proving the above result is motivated by the idea of embedding into $l^1$ an increasing sequence of sub-trees of $\mathcal{T}_w$ chosen to span a sample of $\mu_w$-random vertices of $\mathcal{T}_w$. It is an artifact of the construction of the pair $(\mathcal{T}_w,\mu_w)$ from the excursion $w$ that choosing a $\mu_w$-random sequence of vertices can be related to choosing a collection of uniform random variables from $[0,1]$. Throughout this article, we will use the notation $U=(U_n)_{n\geq 1}$ to represent an independent identically-distributed sequence of $U[0,1]$ random variables built, under $\mathbf{P}$, independently of the random excursion $W$.  We describe fully in Section \ref{crtsec} how a pair of the form $(w,u)\in C([0,1],\mathbb{R}_+)\times[0,1]^\mathbb{N}$, which can be considered to be a particular realisation of $(W,U)$, can be used to construct both $(\mathcal{T}_w,{\mu}_w,{\mathbf{P}}^{\mathcal{T}_w,\mu_w}_\rho)$ and its isometric $l^1$-embedding $(\tilde{\mathcal{T}},\tilde{\mu},\tilde{\mathbf{P}}^\mathcal{T}_\rho)$; at least for a suitably large subset of $C([0,1],\mathbb{R}_+)\times[0,1]^\mathbb{N}$. Similar embeddings are used for discrete trees, see Section \ref{cbptree}.

The reason for choosing ordered trees in Theorem \ref{mainresult} is only for convenience, as it allows us to prove all the convergence results for the finite length sub-trees in an abstract tree space, leaving embedding into $l^1$ until the end, and also, in the annealed result we state below, means we do not have to consider awkward conditional distributions to select these sub-trees. In fact, it is also possible to apply an almost identical argument in the unordered case for any sequence of discrete trees for which the deterministic analogue of the conditions of \cite{Aldous3}, Corollary 19 hold. The stochastic versions of these conditions were used by Aldous in \cite{Aldous3} to demonstrate that it is possible to embed all the relevant objects into $l^1$ in such a way that a random sequence $\{(\mathcal{T}_n,\mu_n)\}_{n\geq 1}$ converges in distribution to $(\mathcal{T}_W,\mu_W)$, the continuum random tree (and associated measure). Thus, as in the ordered case, there are no extra conditions needed to extend from the convergence of trees and measures to the convergence of trees, measures and processes. The only difference in this case is that we will need to use an exchangeability argument similar to \cite{Aldous3}, Theorem 18, to deduce Lemma \ref{subtreelemma}, rather than the excursion one followed here.

After checking the measurability of the embedding $(w,u)\mapsto (\tilde{\mathcal{T}},\tilde{\mu},\tilde{\mathbf{P}}^\mathcal{T}_\rho)$ that we employ as a map from $C([0,1],\mathbb{R}_+)\times[0,1]^\mathbb{N}$ into $\mathcal{K}(l^1)\times\mathcal{M}_1(l^1)\times \mathcal{M}_1(C([0,1],l^1))$, see Section \ref{meassec}, there is no problem in defining a probability law $\mathbb{P}$ on $\mathcal{K}(l^1)\times\mathcal{M}_1(l^1)\times C([0,1],l^1)$ that satisfies
\begin{equation}\label{pdef}
\mathbb{P}\left(A\times B \times C\right)=\int_{C([0,1],\mathbb{R}_+)\times [0,1]^\mathbb{N}} \mathbf{P}((W,U)\in (dw,du))\: \mathbf{1}_{\{\tilde{\mathcal{T}}\in A,\:\tilde{\mu}\in B\}}\tilde{\mathbf{P}}^\mathcal{T}_\rho(C),
\end{equation}
for every measurable $A\subseteq\mathcal{K}(l^1)$, $B\subseteq\mathcal{M}_1(l^1)$, and $C\subseteq C([0,1],l^1)$. In fact, we actually show that it is possible to define a random quintuplet $(W,U,\tilde{\mathcal{T}},\tilde{\mu},\tilde{X})$, where the pair $(\tilde{\mathcal{T}},\tilde{\mu})$ is constructed (measurably) from $(W,U)$, (so that it is simply a random embedding of the continuum random tree and associated measure into $l^1$), in such a way that: the joint law of $(W,U)$ is as described above; the joint law of $(\tilde{\mathcal{T}},\tilde{\mu},\tilde{X})$ is $\mathbb{P}$; and moreover,
\begin{equation}\label{cond}
\mathbf{P}\left(\tilde{X}\in\cdot|\:(W,U)\right)={\mathbf{P}}^{\tilde{\mathcal{T}},\tilde{\mu}}_0,
\end{equation}
where ${\mathbf{P}}^{\tilde{\mathcal{T}},\tilde{\mu}}_0$ is the law of the Brownian motion on $(\tilde{\mathcal{T}},\tilde{\mu})$, started from $0\in\tilde{\mathcal{T}}$. We call the (non-Markovian) process $\tilde{X}$ the Brownian motion on the continuum random tree (isometrically embedded into $l^1$).

A similar law can be constructed in the discrete case. More specifically, let $(\mathcal{T}_n)_{n\geq 1}$ be a sequence of random ordered graph trees with corresponding search-depth functions $(W_n)_{n\geq 1}$, and suppose these are built on our underlying probability space independently of the random variable $U$. Clearly there is a one-to-one correspondence between search-depth functions and $n$-vertex trees, and to imitate the definition of $\mathbb{P}$ we will define the related discrete law in terms of the sequence $(W_n)_{n\geq 1}$. It is straightforward to check that the map from a realisation of a search-depth function and sequence in $[0,1]$, $(w_n,u)$ say, to the $l^1$-embedded triple $(\tilde{\mathcal{T}}_n,\tilde{\mu}_n,\tilde{\mathbf{P}}^{\mathcal{T}_n}_\rho)$ is measurable, and hence we can define a law $\mathbb{P}_n$ on $\mathcal{K}(l^1)\times\mathcal{M}_1(l^1)\times C(\mathbb{R}_+,l^1)$ that satisfies
\begin{equation}\label{pndef}
\mathbb{P}_n\left(A\times B \times C\right)=\int_{C([0,1],\mathbb{R}_+)\times [0,1]^\mathbb{N}} \mathbf{P}((W_n,U)\in (dw,du))\: \mathbf{1}_{\{\tilde{\mathcal{T}}_n\in A,\:\tilde{\mu}_n\in B\}}\tilde{\mathbf{P}}^{\mathcal{T}_n}_\rho(C),
\end{equation}
for every measurable $A\subseteq\mathcal{K}(l^1)$, $B\subseteq\mathcal{M}_1(l^1)$, and $C\subseteq C(\mathbb{R}_+,l^1)$. Similarly to the continuous case, if $(\tilde{\mathcal{T}}_n,\tilde{\mu}_n,\tilde{X}^n)$ represents a random variable with law $\mathbb{P}_n$, then $(\tilde{\mathcal{T}}_n,\tilde{\mu}_n)$ is equal in distribution to a certain random $l^1$-embedding of $(\mathcal{T}_n,\mu_n)$. Moreover, conditional on $\tilde{\mathcal{T}}_n$, the process $\tilde{X}^n$ is a simple random walk on the elements of $\tilde{\mathcal{T}}_n$ (edges are assumed to be between points separated by a unit distance) started from the origin.

We are now ready to state our annealed convergence result. The rescaling operator $\Theta_n$ is redefined on $\mathcal{K}(l^1)\times \mathcal{M}_1(l^1)\times C(\mathbb{R}_+,l^{1})$ in the obvious way, so that if $(\tilde{K},\tilde{\nu},\tilde{f})$ is an element of this space, then $\Theta_n(\tilde{K},\tilde{\nu},\tilde{f}):=(n^{-1/2}\tilde{K},\tilde{\nu}(n^{1/2}\cdot),(n^{-1/2}\tilde{f}(tn^{3/2}))_{t\in[0,1]})$.
The notation $\Rightarrow$ is used in all that follows to represent convergence in distribution.

{\thm\label{ann} Suppose that $(\mathcal{T}_n)_{n\geq 1}$ is a sequence of random ordered graph trees whose search-depth functions $(W_n)_{n\geq 1}$ satisfy
\begin{equation}\label{sdc}
n^{-1/2}W_n\Rightarrow W
\end{equation}
in $C([0,1],\mathbb{R}_+)$, then if $(\mathbb{P}_n)_{n\geq 1}$ and $\mathbb{P}$ are probability measures satisfying (\ref{pdef}) and (\ref{pndef}), respectively, then
\[\mathbb{P}_n\circ\Theta_n^{-1}\rightarrow \mathbb{P}\]
weakly as measures on the space $\mathcal{K}(l^1)\times\mathcal{M}_1(l^1)\times C([0,1],l^1)$.}
\bigskip

Equivalently, we can also write this result in terms of random variables.

{\cor Assume that, for each $n$, the law of the random triple $(\tilde{\mathcal{T}}_n,\tilde{\mu}_n,\tilde{X}^n)$, which consists of a random $l^1$-embedded graph tree, measure and associated simple random walk, is given by $\mathbb{P}_n$, and $(\tilde{\mathcal{T}},\tilde{\mu},\tilde{X})$ is the random embedding into $l^1$ of the continuum random tree and Brownian motion upon it (so that it has law $\mathbb{P}$). If $n^{-1/2}W_n\Rightarrow W$, we have that
\[\Theta_n\left(\tilde{\mathcal{T}}_n,\tilde{\mu}_n,\tilde{X}^n\right)\Rightarrow \left(\tilde{\mathcal{T}},\tilde{\mu},\tilde{X}\right)\]
in the space $\mathcal{K}(l^1)\times\mathcal{M}_1(l^1)\times C([0,1],l^1)$.}
\bigskip

As a final remark, we note that there is nothing particularly special or fundamental about the space $l^1$, and there should be no problem in stating the results of this article in a more abstract space of metric space trees, measures and processes on them by, for example, generalising the spaces investigated in \cite{EPW} and \cite{Greven}. Due to the length of the article, we leave such a presentation for future work.

This article is almost entirely devoted to demonstrating Theorem \ref{mainresult}. After introducing the majority of the notation we use and some relevant background material in Section \ref{prelim}, we provide an overview of the proof in Section \ref{overv}, which explains how the argument is structured. In Section \ref{meassec} we tackle various measurability issues, and the results we prove there allow us to derive from Theorem \ref{mainresult} the remaining conclusions of this section.

\section{Preliminaries}\label{prelim}

\subsection{Abstract trees and projections}\label{abstracttree}

Although the conclusion of this article is stated in terms of trees embedded into $l^1$, for most of the arguments we do not need to be this specific about the space in which we are working. In this section, we introduce some notation and concepts for arbitrary metric space trees. For the purposes of this and the next section, we shall assume that $K$ is a dendrite, which means that it is an arc-wise connected topological space, containing no subset homeomorphic to the circle. We shall also suppose that $d_K$ is a shortest path metric on $K$, which means that it is additive along the (non-self intersecting) paths of $K$. In this article, we shall have cause to refer to the root of various graph trees/dendrites. This is a distinguished vertex, and we shall denote it by $\rho$. For brevity, we will usually write the triple $(K,d_K,\rho)$ as simply $K$. A metric space of this form is also known as a rooted real tree. Note that much of the notation and terminology we introduce here for the dendrite $K$ also makes sense for graph trees, and so we shall apply it to graphs with no further explanation.

One of the consequences of $K$ being a dendrite is that, for any $x,y\in K$, there exists a unique (non-self intersecting) path from $x$ to $y$. We will use the notation $[[x,y]]$ to represent such a path. Furthermore, between any 3 vertices $x,y,z\in K$ there is a unique branch point $b^K(x,y,z)\in K$ which satisfies
\[\{b^K(x,y,z)\}=[[x,y]]\cap[[y,z]]\cap[[z,x]].\]
We define the degree of a vertex $x\in K$ by
\[\mathrm{deg}_K(x):=\#\{\mbox{connected components of }{K}\backslash\{x\}\},\]
which takes values in $\mathbb{N}\cup\{0,\infty\}$.

In the analysis of the stochastic processes that follows in later sections, we will use the idea of observing processes on reduced sub-trees (strictly speaking, these are reduced sub-dendrites). Given $A\subseteq K$, the reduced sub-tree $r(K,A)$ is the smallest path-wise connected subset of $K$ containing $A\cup\{\rho\}$. In particular, we have
\[r(K,A):=\bigcup_{x\in A} [[\rho,x]].\]
The subset $r(K,A)$ is clearly a dendrite, and in the case of $A$ being finite, $r(K,A)$ is a closed subset of $K$.

Given an arbitrary closed sub-tree of $K$, that is a closed set $K'\subseteq K$ such that $(K',d_K)$ is a dendrite, there is a natural projection from $K$ onto $K'$. This continuous map will be denoted by $\phi_{K,K'}$, and may be defined in the following way: for a point $x\in K$, $\phi_{K,K'}(x)$ is the unique point in $[[\rho,x]]$ such that
\begin{equation}\label{project}
[[\phi_{K,K'}(x),x]]\cap K'=\{\phi_{K,K'}(x)\}.
\end{equation}
Note that, necessarily, $\phi_{K,K'}(x)\in K'$. Perhaps a clearer way of describing the projection is provided by the observation that, for $x\in K$, $\phi_{K,K'}(x)$ is the point in $K'$ closest to $x$.

We now provide a brief introduction to the connection between trees and excursions. This is an area which has been of much recent interest and we shall use the idea to define the continuum random tree in Section \ref{crtsec}. First, let $\mathcal{W}$ be the collection of continuous functions $w:\mathbb{R}_+\rightarrow \mathbb{R}_+$ for which there exists a $\tau(w)>0$ such that $w(t)>0$ if and only if $t\in (0,\tau(w))$. The set $\mathcal{W}$ is the space of excursions. For future use, we introduce the notation $\mathcal{W}^{(1)}:=\{w\in \mathcal{W}:\:\tau(w)=1\}$ to represent the excursions of length 1. Given a function $w\in \mathcal{W}$, we define a distance on $[0,\tau(w)]$ by setting
\[d_w(s,t):=w(s)+w(t)-2m_w(s,t),\]
where $m_w(s,t):=\inf\{w(r):\:r\in[s\wedge t,s\vee t]\}$. Then, we use the equivalence,
\begin{equation}\label{equivrel}
s\sim t\hspace{20pt}\Leftrightarrow\hspace{20pt}d_w(s,t)=0
\end{equation}
to define $\mathcal{T}_w:=[0,\tau(w)]/\sim$. We can write this as $\mathcal{T}_w=\{[s]:\:s\in[0,\tau(w)]\}$, where $[s]$ is the equivalence class containing $s$. It is then elementary (see \cite{LegallDuquesne}, Section 2.1) to check that $d_{\mathcal{T}_w}([s],[t]):=d_w(s,t)$, defines a metric on $\mathcal{T}_w$, and also that $\mathcal{T}_w$ is a compact dendrite. Furthermore, the metric $d_{\mathcal{T}_w}$ is a shortest path metric on $\mathcal{T}_w$. The root of the tree $\mathcal{T}_w$ is defined to be the equivalence class $[0]$.

A natural measure to impose upon $\mathcal{T}_w$ is the projection of Lebesgue measure on $[0,\tau(w)]$. For open $A\subseteq\mathcal{T}_w$, let
\begin{equation}\label{muw}
\mu_w(A):=\lambda\left(\{t\in[0,\tau(w)]:\:[t]\in A\}\right),
\end{equation}
where, throughout this article, $\lambda$ is the usual 1-dimensional Lebesgue measure. This defines a Borel measure on $(\mathcal{T}_w,d_{\mathcal{T}_w})$, with total mass equal to $\tau(w)$.

To complete this section, we explain how to use a sequences $u=(u_n)_{n\geq 1}\in[0,\tau(w)]^\mathbb{N}$ to define a sequence of increasing sub-trees of $\mathcal{T}_w$. First, define a collection of vertices of $\mathcal{T}_w$ by
\begin{equation}\label{vertices}
v_n:=[u_n],
\end{equation}
where $[t]$ is the equivalence class of $t\in[0,\tau(w)]$, as defined above. From this collection of vertices we obtain a sequence of closed sub-trees of $\mathcal{T}_w$ by defining, for $k\geq 1$,
\[\mathcal{T}_{w,u}(k):=r(\mathcal{T}_w,\{v_1,\dots,v_k\}).\]
Note that this sequence is increasing in the sense that $\mathcal{T}_{w,u}(k)\subseteq\mathcal{T}_{w,u}(k+1)$, for every $k$. The projection of $\mu_w$ onto $\mathcal{T}_{w,u}(k)$ will be denoted
\[\mu_{w,u}^{(k)}:=\mu\circ\phi_{\mathcal{T}_w,\mathcal{T}_{w,u}(k)}^{-1}.\]
This will not be the only measure of interest on $\mathcal{T}_{w,u}(k)$. Since $\mathcal{T}_{w,u}(k)$ is a tree consisting of a finite number of edges with strictly positive total edge length, there is no problem in defining Lebesgue measure $\lambda_{w,u}^{(k)}$ on $\mathcal{T}_{w,u}(k)$. More specifically, this is the measure that satisfies
\begin{equation}\label{leb}
\lambda_{w,u}^{(k)}([[x,y]])\propto d_{\mathcal{T}_w}(x,y),\hspace{20pt}\forall x,y,\in\mathcal{T}_{w,u}(k).
\end{equation}
We shall normalise $\lambda_{w,u}^{(k)}$ so that it is a probability measure on $\mathcal{T}_{w,u}(k)$. We remark that this is indeed possible by applying the fact that $\mathrm{diam}\mathcal{T}_w$ is finite, which is a simple consequence of the compactness of $\mathcal{T}_w$. Both $\mu_{w,u}^{(k)}$ and $\lambda_{w,u}^{(k)}$ are clearly Borel measures on $\mathcal{T}_{w,u}(k)$, and it is straightforward to check that $\mu_{w,u}^{(k)}(A)>0$ and $\lambda_{w,u}^{(k)}(A)>0$ for every non-empty open $A\subseteq\mathcal{T}_w(k)$.

Note that we will usually drop the subscripts $w$ and $u$ from the objects described above when it is clear which excursion and sequence is being considered.

\subsection{Processes on abstract trees}\label{abstractproc}

Using a result of Kigami, it is possible to establish the existence of ``nice'' Markov processes on a wide class of dendrites. As in the previous section, we assume that $(K,d_K)$ is a dendrite equipped with a shortest path metric. We shall suppose further that $\nu$ is a $\sigma$-finite Borel measure on $K$ that satisfies $\nu(A)>0$ for every non-empty open set $A\subseteq K$. The following result is proved by Kigami as Theorem 5.4 of \cite{Kigamidendrite}. Definition 0.5 of \cite{Kigamidendrite} specifies the precise conditions that make a symmetric, non-negative quadratic form a finite resistance form. For more examples of this type of form, see \cite{Kigami}. We shall not explain here how to construct the finite resistance form associated with a shortest path metric on a dendrite, as knowledge of this is non-essential for the results of this article. Full details are given in Section 3 of \cite{Kigamidendrite}. We shall however, continue to use the notation $(\mathcal{E}_K,\mathcal{F}_K)$ to represent such a form.

{\lem Suppose $(K,d_K)$ is locally compact and complete, then $(\mathcal{E}_K,\mathcal{F}_K\cap L^2(K,\nu))$, where $(\mathcal{E}_K,\mathcal{F}_K)$ is the finite resistance form associated with $(K,d_K)$, is a local, regular Dirichlet form on $L^2(K,\nu)$.}
\bigskip

We now describe briefly the natural construction of the Markov process corresponding to $(\mathcal{E}_K, \mathcal{F}_K)$ and measure $\nu$, and outline the properties of this process that will be relevant to this article. Given the Dirichlet form $(\frac{1}{2}\mathcal{E}_K,\mathcal{F}_K\cap L^2(K,\nu))$, we can use the standard association to define a non-negative self-adjoint operator, $-\Delta_K$, which has domain dense in $L^2(K,\nu)$ and satisfies
\[\mbox{$\frac{1}{2}$}\mathcal{E}_K(u,v)=-\int_K u\Delta_{K} vd\nu,\hspace{20pt}\forall u\in\mathcal{F}_K\cap L^2(K,\nu),v\in\mathcal{D}(\Delta_K).\]
Although the factor $\frac{1}{2}$ looks rather awkward here, it will be useful in ensuring a particular time-scaling for the reversible Markov process,
\[X^{K,\nu}=((X^{K,\nu}_t)_{t\geq 0}, \mathbf{P}_x^{K,\nu}, x\in K),\]
which is defined from the semi-group given by $P_t:=e^{t\Delta_K}$. In fact, the locality of our Dirichlet form ensures that the process $X^{K,\nu}$ is a diffusion on $K$.

In the case when $(K,d_K)$ is compact, Aldous defines in \cite{Aldous2} a Brownian motion on $(K,d_K,\nu)$ to be a process with the following properties. Note first that, since we only use one metric on any particular dendrite, we will omit the metric from the notation from now on.
\newcounter{listcount}
\begin{list}{\roman{listcount})}{\usecounter{listcount} \setlength{\rightmargin}{\leftmargin}}
\item Continuous sample paths.
\item Strong Markov.
\item Reversible with respect to its invariant measure $\nu$.
\item For $x,y\in K$, $x\neq y$, we have
\[\mathbf{P}_{z}^{K,\nu}\left( \sigma_{x}<\sigma_y \right)=\frac{d_K(b^K(z,x,y),y)}{d_K(x,y)},\hspace{20pt}\forall z\in K,\]
where $\sigma_{x'}:=\inf\{t>0:\:X^{K,\nu}_t=x'\}$ is the hitting time of $x'\in K$.
\item For $x,y\in K$, the mean occupation measure for the process started at $x$ and killed on hitting $y$ has density
\[2d_K(b^K(z,x,y),y)\nu(dz),\hspace{20pt}\forall z\in K.\]
\end{list}
As remarked in Section 5.2 of \cite{Aldous2}, these properties are enough to guarantee the uniqueness of Brownian motion on $(K,\nu)$. We now discuss existence. In fact, the following proposition was essentially proved in \cite{Croydoncrt} and gives us that the process constructed from the Dirichlet form associated with $(\mathcal{E}_K,\mathcal{F}_K)$ and $\nu$, as above, is actually the Brownian motion on $(K,\nu)$. Note how, in this result, the domain of the Dirichlet form does not depend on the choice of measure. Since it can be proved using exactly the same arguments as Section 8 of \cite{Croydoncrt}, we simply state the result.

{\propn \label{BM} Let $(K,d_K)$ be a compact dendrite and $(\mathcal{E}_K,\mathcal{F}_K)$ be the finite resistance form associated with $(K,d_K)$. Then $(\frac{1}{2}\mathcal{E}_K,\mathcal{F}_K)$ is a local, regular Dirichlet form on $L^2(K,\nu)$, and furthermore, the corresponding Markov process $X^{K,\nu}$ is Brownian motion on $(K,\nu)$.}

\subsection{Continuum random tree properties}\label{crtsec}

In this section, we introduce the continuum random tree and a certain collection of random sub-trees of it. Our starting point is that we assume that we are given a pair of random variables $(W,U)$ built on an underlying probability space with probability measure $\mathbf{P}$. Under $\mathbf{P}$, the process $W=(W_t)_{t\in [0,1]}$ is a normalised Brownian excursion. For a precise description of the law of $W$, see \cite{RevuzYor}, Chapter XII. The random variable $U=(U_n)_{n\geq 1}$ is a sequence of independent $U[0,1]$ random variables, independent of $W$.

Since the random variable $(W,U)$ takes values in $\mathcal{W}^{(1)}\times [0,1]^{\mathbb{N}}$, $\mathbf{P}$-a.s., we can use the procedure of Section \ref{abstracttree} to define (on at least on a set of probability 1) a compact dendrite and measure, $(\mathcal{T}_W,d_{\mathcal{T}_W},\mu_W)$, an increasing sequence of sub-trees $(\mathcal{T}_{W,U}(k))_{k\geq 1}$, and also, for each $k\geq 1$, the measures $\mu_{W,U}^{(k)}$ and $\lambda_{W,U}^{(k)}$. The dendrite $\mathcal{T}_W$ is the continuum random tree. We shall, in future, drop the subscripts $W$ and $U$ when it will not cause confusion. We note that $\tau(W)=1$, $\mathbf{P}$-a.s., and so $\mu$ is a probability measure on $\mathcal{T}$, $\mathbf{P}$-a.s.

In analysing random variables which take values in infinite dimensional spaces, such as continuous time stochastic processes, it is often useful to proceed by investigating finite dimensional distributions and taking some limit. As discussed in Section 2.4 of \cite{Aldous2}, the natural substitute for this in proving the convergence of discrete trees to the continuum random tree is the investigation of random finite dimensional distributions. We briefly note that this is the technique that we apply, since if we denote the sequence of vertices that are used to construct $\mathcal{T}(k)$ from $\mathcal{T}$ by $V_n:=[U_n]$, then conditional on $W$, by definition of the measure $\mu$ as the projection of Lebesgue measure on $[0,1]$ onto $\mathcal{T}$, we have that the vertices $(V_n)_{n\geq 1}$ are an independent, identically distributed sample of $\mu$-distributed random variables.

In the following lemma, we collect together some important properties of $(\mathcal{T},\mu)$ and the sequence of sub-trees and measures.

{\lem \label{crtprops} There exists a measurable set $\Gamma\subseteq \mathcal{W}^{(1)}\times[0,1]^{\mathbb{N}}$ such that $\mathbf{P}((W,U)\in\Gamma)=1$, and also if $(w,u)\in\Gamma$, then\\
(a) $\mathcal{T}$ is a compact dendrite with shortest path metric $d_\mathcal{T}$.\\
(b) For every $x\in\mathcal{T}$, $\mathrm{deg}_\mathcal{T}(x)\leq 3$.\\
(c) If $N(\mathcal{T},\varepsilon)$ is the number of $\varepsilon$ balls needed to cover $\mathcal{T}$, then
\begin{equation}\label{coversize}
\limsup_{\varepsilon\rightarrow 0}\varepsilon^2 N(\mathcal{T},\varepsilon)<\infty.
\end{equation}
(d) The elements of $(u_n)_{n\geq 1}$ are disjoint, and so are the elements of the collection $(v_n)_{n\geq 1}$, as defined at (\ref{vertices}). Moreover, the collection of vertices $(v_n)_{n\geq 1}$ is dense in $\mathcal{T}$.\\
(e) The Brownian motion on $(\mathcal{T},\mu)$ exists, and admits a heat kernel $(p_t(x,y))_{t>0,x,y\in\mathcal{T}}$ that satisfies
\begin{equation}\label{tdupper}
\limsup_{t\rightarrow0}t^{2/3}(\ln t^{-1})^{-1/3}\sup_{x\in\mathcal{T}}p_t(x,x)<\infty.
\end{equation}
(f) For each $k$, $\mu^{(k)}$ and $\lambda^{(k)}$ are Borel probability measures on $\mathcal{T}(k)$ that satisfy $\mu^{(k)}(A)>0$ and $\lambda^{(k)}(A)>0$ for every non-empty open $A\subseteq\mathcal{T}(k)$.\\
(g) As $k\rightarrow \infty$, $\lambda^{(k)}\rightarrow \mu$ weakly as Borel probability measures on $\mathcal{T}$.}
\begin{proof} Parts (a) and (f) are obvious from the construction in Section \ref{abstracttree}. Parts (b) and (c) are covered by \cite{LegallDuquesne}, Theorem 4.6(iv) and Proposition 5.2 respectively. The proof of part (d) requires only elementary analysis, and is therefore omitted. The existence of a heat kernel for Brownian motion on $(\mathcal{T},\mu)$ was established in \cite{Croydoncrt}; the estimate of part (e) was also proved in the same reference. For part (g), see \cite{Aldous1}, Theorem 3(ii). \end{proof}

In future, we shall fix a particular $\Gamma\subseteq \mathcal{W}^{(1)}\times[0,1]^{\mathbb{N}}$ that satisfies the claims of the above lemma. We shall denote by
\begin{equation}\label{wstar}
\mathcal{W}^*:=\left\{w:\:(w,u)\in\Gamma\mbox{ for some }u\in[0,1]^{\mathbb{N}}\right\}
\end{equation}
the projection onto the first coordinate of $\Gamma$. Roughly speaking, this represents a set of typical realisations of the continuum random tree $(\mathcal{T},\mu)$ that can be approximated in a good way by a collection $(\mathcal{T}(k),\lambda^{(k)})$ or $(\mathcal{T}(k),\mu^{(k)})$ of suitably selected sub-trees. Clearly $\mathbf{P}(W\in \mathcal{W}^*)$=1.

Before continuing, we derive an extra tightness condition that holds when $(\mathcal{T}(k))_{k\geq 1}$ and $\mathcal{T}$ are constructed from $(w,u)\in\Gamma$.

{\lem \label{deltakdecay} For $(w,u)\in\Gamma$, we have that, as $k\rightarrow\infty$,
\[\Delta^{(k)}:=\sup_{x\in\mathcal{T}}d_\mathcal{T} (x,\phi_{\mathcal{T},\mathcal{T}(k)}(x))\rightarrow 0.\]}
\begin{proof} Fix $\varepsilon>0$. By the compactness of $\mathcal{T}$, there exists a finite collection, $(x_i)_{i=1}^N$, of elements of $\mathcal{T}$ such that $\mathcal{T}\subseteq \bigcup_{i=1}^N B(x_i, \varepsilon/2)$. Furthermore, by the denseness of $(v_k)_{k\geq 1}$, for each $x_i$, we can find a $k_i$ such that $d_\mathcal{T}(x_i, v_{k_i})\leq\varepsilon/2$. Now suppose $k\geq k_0:=\max_{i\in\{1,\dots,N\}}k_i$ and $x\in \mathcal{T}$. Since by definition, $\mathcal{T}(k_0)\subseteq\mathcal{T}(k)$ and $\phi_{\mathcal{T},\mathcal{T}(k)}(x)$ is the point of $\mathcal{T}(k)$ closest to $x$, we have $d_\mathcal{T}(x,\phi_{\mathcal{T},\mathcal{T}(k)}(x))\leq d_\mathcal{T}(x,\phi_{\mathcal{T},\mathcal{T}(k_0)}(x)).$ Also, by choice of $(x_i)_{i=1}^N$, we must have, $x\in B(x_i,\varepsilon/2)$ for some $i$. Applying this, and using the fact that the branch point $b^\mathcal{T}(0,x_i,v_{k_i})$ is necessarily an element of $\mathcal{T}(k_0)$, it may be deduced from the previous inequality that
\begin{eqnarray*}
d_\mathcal{T}(x,\phi_{\mathcal{T},\mathcal{T}(k)}(x))&\leq& d_\mathcal{T}(x,b^\mathcal{T}(0,x_i,v_{k_i}))\\
&\leq&\frac{\varepsilon}{2}+d_\mathcal{T}(x_i,b^\mathcal{T}(0,x_i,v_{k_i}))\\
&\leq& \frac{\varepsilon}{2}+d_\mathcal{T}(x_i,v_{k_i})\\
&\leq&\varepsilon.
\end{eqnarray*}
Thus $\Delta^{(k)}\leq \varepsilon$, for all $k\geq k_0$, and the lemma follows.
\end{proof}

We now explain how to embed $(\mathcal{T},\mu)$ into $l^1$ by using the sequence of sub-trees $(\mathcal{T}(k))_{k\geq 1}$. For the purposes of the following discussion, we assume that $(w,u)\in \Gamma$. For each $k\in\mathbb{N}$, define $\tilde{\mathcal{T}}(k)$ to be the subset of $l^1$ obtained by isometrically embedding $\mathcal{T}(k)$ into $l^1$ from the sequence of vertices $(v_1,\dots,v_k)$ using the sequential construction of \cite{Aldous3}, Section 2.2. In short, this procedure involves adding successive branches orthogonally. More precisely, set $\tilde{\mathcal{T}}(1):=\{tz_1:\:t\in[0,d_\mathcal{T}(\rho,v_1)]\}$, where $(z_k)_{k\geq 1}$ is the canonical basis for $l^1$. Suppose all the sets $(\tilde{\mathcal{T}}(k'))_{k'\leq k}$ are defined, then there exists isometries $\psi^{(k')}:\mathcal{T}(k')\rightarrow\tilde{\mathcal{T}}(k')$, for $k'\leq k$, which may be uniquely determined by insisting that $\psi^{(k')}\vline_{\mathcal{T}(k'-1)}=\psi^{(k'-1)}$ for each $k'\leq k$, where we use the convention that $\mathcal{T}(0)=\{\rho\}$ and $\tilde{\mathcal{T}}(0)=\{0\}$. The inductive step is the following:
\[\tilde{\mathcal{T}}(k+1):=\tilde{\mathcal{T}}(k)\cup\left(\psi^{(k)}(\phi_{\mathcal{T},\mathcal{T}(k)}(v_{k+1}))+\{tz_{k+1}:\:t\in[0,d_\mathcal{T}(\phi_{\mathcal{T},\mathcal{T}(k)}(v_{k+1}),v_{k+1})]\}\right).\]
It is easy to check that this procedure results in an increasing sequence of subsets of $l^1$ such that, for each $k$, $(\tilde{\mathcal{T}}(k),\|\cdot-\cdot\|)$ is an isometric copy of $(\mathcal{T}(k),d_\mathcal{T})$. We will denote $\tilde{\mu}^{(k)}:=\mu^{(k)}\circ{\psi^{(k)}}^{-1}$ and $\tilde{\lambda}^{(k)}:=\lambda^{(k)}\circ{{\psi^{(k)}}^{-1}}$, which are Borel probability measures on $l^1$.

Since on $\Gamma$ the vertices $(v_k)_{k\geq 1}$ are dense in $\mathcal{T}$, it is a simple exercise to define a distance-preserving map $\psi:\mathcal{T}\rightarrow l^1$ that satisfies $\psi\vline_{\mathcal{T}(k)}=\psi^{(k)}$ for each $k$. Denoting $\tilde{\mathcal{T}}:=\psi(\mathcal{T})$, we have that $(\tilde{\mathcal{T}},\|\cdot-\cdot\|)$ is an isometric copy of $(\mathcal{T},d_\mathcal{T})$. Moreover, if we define the root of $\tilde{\mathcal{T}}$ to be 0, then $\psi$ is root-preserving. Also define $\tilde{\mu}:=\mu\circ\psi^{-1}$ and $\tilde{\mathbf{P}}^\mathcal{T}_\rho:={\mathbf{P}}^{\mathcal{T},\mu}_\rho\circ \psi^{-1}$, where ${\mathbf{P}}^{\mathcal{T},\mu}_\rho$ is the law of the Brownian motion on $(\mathcal{T},\mu)$ started from the root. Finally, although there is no problem with defining the objects $\tilde{\mathcal{T}}$, $\tilde{\mu}$, $\tilde{\mathbf{P}}^\mathcal{T}_\rho$, $\tilde{\mathcal{T}}(k)$, $\tilde{\mu}^{(k)}$ and  $\tilde{\lambda}^{(k)}$ in a deterministic way for each $(w,u)\in\Gamma$, to prove the distributional result of Theorem \ref{ann} and the conditional relation at (\ref{cond}) we need to show that the construction is $(W,U)$-measurable, and we do this in Section \ref{meassec}.

\subsection{Coupling processes on the continuum random tree}\label{crtproc}

In this section we suppose that $(w,u)\in\Gamma$ is fixed. On $\Gamma$, the assumptions of Proposition \ref{BM} hold for both $(\mathcal{T},\mu)$ and $(\mathcal{T}(k),\lambda^{(k)})$, and so we can construct the Brownian motions on these spaces. It will be useful to couple these processes, and so we will construct the Brownian motions on $(\mathcal{T}(k),\lambda^{(k)})$ using a simple time-change argument.

First, denote by $X$ the Brownian motion on $(\mathcal{T},\mu)$ started from $\rho$, so that, under $\mathbf{P}$, the law of $X$ is $\mathbf{P}^{\mathcal{T},\mu}_\rho$. We start by showing that this process admits $\mathbf{P}$-a.s. jointly continuous local times. The argument we use follows closely that of \cite{Barlow}, Theorem 7.21, in which the corresponding result was proved for the diffusions on certain deterministic post-critically finite, self-similar fractals.

{\lem \label{localtimecont} For $(w,u)\in\Gamma$, there exist local times $(L_t(x))_{t\geq 0, x\in\mathcal{T}}$ for the process $(X_t)_{t\geq 0}$ which are $\mathbf{P}$-a.s. jointly continuous in $t$ and $x$.}
\begin{proof} The existence of jointly measurable local times for Brownian motion on $(\mathcal{T},\mu)$ was essentially demonstrated in the proof of \cite{Croydoncrt}, Lemma 8.2, and so it remains to show continuity. Recall that on $\Gamma$, the Brownian motion on $(\mathcal{T},\mu)$ admits a transition density that satisfies the upper estimate at (\ref{tdupper}). The 1-potential density, $u$, is defined from the transition density by $u(x,y):=\int_0^{\infty}e^{-t}p_t(x,y)dt$, and it is easily deduced from (\ref{tdupper}) that $u(x,y)$ is finite for all $x,y\in\mathcal{T}$. This allows us to apply \cite{MarcusRosen}, Theorem 1, to deduce that the $\mathbf{P}$-a.s. continuity of the local times of $X$ is equivalent to the $\mathbf{P}$-a.s. continuity of the process $(G(x))_{x\in\mathcal{T}}$, which is defined to be a mean zero, Gaussian process with covariances given by $(u(x,y))_{x,y\in\mathcal{T}}$. However, to prove the continuity of $(G(x))_{x\in\mathcal{T}}$, by applying \cite{DudleyGauss}, Theorem 2.1, it is sufficient to show that $\int_{0}^1\sqrt{\ln N(\mathcal{T},\varepsilon)}d\varepsilon<\infty$, where $N(\mathcal{T},\varepsilon)$ is the smallest number of balls of radius $\varepsilon$ needed to cover $\mathcal{T}$. On $\Gamma$, we have by (\ref{coversize}) that the relevant integral is indeed finite, and so the proof is complete.
\end{proof}

We now explain the coupling that we will apply. For $k\geq 1$, define the continuous additive functional $(A^{(k)}_t)_{t\geq 0}$ by
\begin{equation}\label{akdef}
A^{(k)}_t:=\int_{{\mathcal{T}}(k)}L_t(x){\lambda}^{(k)}(dx)
\end{equation}
and its inverse by
\begin{equation}\label{taukdef}
\tau^{(k)}(t):=\inf\{s:\:A^{(k)}_s> t\}.
\end{equation}
The process $(B^{(k)}_t)_{t\geq 0}$ is then defined by setting
\begin{equation}\label{bkdef}
B^{(k)}_t:=X_{{\tau}^{(k)}(t)}.
\end{equation}
In the following lemma, we use the trace theorem for Dirichlet forms to deduce that, under $\mathbf{P}$, $B^{(k)}$ is the Brownian motion on $({\mathcal{T}}(k),{\lambda}^{(k)})$ started from $\rho$.

{\lem Fix $(w,u)\in\Gamma$ and $k\in\mathbb{N}$. Under $\mathbf{P}$, the process $B^{(k)}$ has law $\mathbf{P}_\rho^{\mathcal{T}(k),\lambda^{(k)}}$.}
\begin{proof} Fix $k\geq 1$, and let $(\mathcal{E}_{{\lambda}^{(k)}},\mathcal{F}_{{\lambda}^{(k)}})$ be the trace of $(\mathcal{E}_\mathcal{T},\mathcal{F}_\mathcal{T})$ onto $\mathcal{T}(k)$ with respect to the measure ${\lambda}^{(k)}$, where $(\mathcal{E}_\mathcal{T},\mathcal{F}_\mathcal{T})$ is the finite resistance form associated with $(\mathcal{T},d_\mathcal{T})$. In particular, we set $\mathcal{E}_{\lambda^{(k)}}(u,u):=\inf\{\mathcal{E}_\mathcal{T}(v,v):\:v|_{\mathcal{T}(k)}=u,\:v\in \mathcal{F}_\mathcal{T}\}$ for $u\in L^2(\mathcal{T}(k),\lambda^{(k)})$, and let $\mathcal{F}_{\lambda^{(k)}}$ be the set of functions for which this infimum exists finitely. By the trace theorem for Dirichlet forms, see \cite{FOT}, Theorem 6.2.1., we have that $B^{(k)}$ is the Markov process associated with $(\frac{1}{2}\mathcal{E}_{\lambda^{(k)}},\mathcal{F}_{\lambda^{(k)}})$, considered as a Dirichlet form on $L^2(\mathcal{T}(k),\lambda^{(k)})$, started from $\rho$.

Since $(\mathcal{E}_\mathcal{T},\mathcal{F}_\mathcal{T})$ is a finite resistance form, it is straightforward to check that so is $(\mathcal{E}_{\lambda^{(k)}},\mathcal{F}_{\lambda^{(k)}})$. Hence, because a finite resistance form is determined by its effective resistance metric (see \cite{Kigamidendrite}, Definition 0.5, for a precise definition of such a metric, and Section 3 for the correspondence), and we can easily check that the relevant metric is simply $d_\mathcal{T}$ restricted to $\mathcal{T}(k)$, we have
$(\mathcal{E}_{\lambda^{(k)}},\mathcal{F}_{\lambda^{(k)}})=(\mathcal{E}_{\mathcal{T}(k)},\mathcal{F}_{\mathcal{T}(k)})$,
where $(\mathcal{E}_{\mathcal{T}(k)},\mathcal{F}_{\mathcal{T}(k)})$ is the finite resistance form associated with $(\mathcal{T}(k),d_\mathcal{T})$. Thus $B^{(k)}$ is the Markov process associated with $(\frac{1}{2}\mathcal{E}_{\mathcal{T}(k)},\mathcal{F}_{\mathcal{T}(k)})$, considered as a Dirichlet form on $L^2(\mathcal{T}(k),\lambda^{(k)})$, started from $\rho$. Consequently Proposition \ref{BM} implies that it is Brownian motion on $(\mathcal{T}(k),\lambda^{(k)})$ started from $\rho$, as claimed.
\end{proof}

\subsection{Discrete trees}\label{cbptree}

In this section, we shall describe the notation that we will use for discrete trees. Since the excursion description of discrete trees is well documented in \cite{Aldous3}, we shall not present the full details, but simply highlight the results which will be important here. First, let $(\mathcal{T}_n)_{n\geq 1}$ be a collection of (rooted) ordered graph trees on $n$ vertices, and, for each $n$, define the function $\hat{w}_n:\{1,\dots,2n-1\}\rightarrow\mathcal{T}_n$ to be the depth-first search around $\mathcal{T}_n$. We extend $\hat{w}_n$ so that $\hat{w}_n(0)=\hat{w}_n(2n)=\rho$, where $\rho=\rho(\mathcal{T}_n)$ is the root of $\mathcal{T}_n$. Define the search-depth process $w_n$ by
\[w_n(i/2n):=d_{\mathcal{T}_n}(\rho,\hat{w}_n(i)),\hspace{20pt}0\leq i\leq 2n,\]
where $d_{\mathcal{T}_n}$ is the graph distance on $\mathcal{T}_n$. Also, extend the definition of $w_n$ to the whole of the interval $[0,1]$ by linear interpolation, so that $w_n$ takes values in $C([0,1],\mathbb{R}_+)$.

Analogously to the definition of $\mathcal{T}(k)$ from $(w,u)$, we shall construct a collection of trees which are sub-trees of $\mathcal{T}_n$ spanning $k$ vertices (some possibly repeated) from $w_n$ and a sequence in $[0,1]^\mathbb{N}$. For each $n$, we denote by $u^n=(u^n_k)_{k\geq 1}$ an element of $[0,1]^\mathbb{N}$. We will presuppose the following condition throughout the remainder of this section, and frequently in subsequent sections, which formalises the notion that $(n^{-1/2}w_n,u^n)$ converges to a typical realisation of the random variable $(W,U)$.

\begin{ass}{\it For each $n$, the sequence $(u^n_k)_{k\geq 1}$ is dense in $[0,1]$, and also
\[\left(n^{-1/2}w_n,u^n\right)\rightarrow (w,u),\]
in $C([0,1],\mathbb{R}_+)\times [0,1]^\mathbb{N}$, for some $(w, u)\in \Gamma$.}
\end{ass}

Since $\hat{w}_n$ is only defined at integer values, we require a slightly more complicated procedure to allow us to use $u^n$ to choose from the $n$ vertices of $\mathcal{T}_n$. For each $n\geq 0$, define the function $\gamma_n:[0,1]\rightarrow[0,1]$ by setting
\[\gamma_n(t):=\left\{\begin{array}{ll}
                        \lfloor 2nt \rfloor/2n, & \mbox{if }w_n(\lfloor 2nt \rfloor/2n)\geq w_n( \lceil 2nt \rceil/2n),\\
                        \lceil 2nt \rceil/2n, & \mbox{otherwise.}
                      \end{array}
      \right.\]
The reason for introducing this particular function is that, by applying an argument similar to Lemma 12 of \cite{Aldous3}, it is possible to show that if $U_1$ is uniform on $[0,1]$, with respect to Lebesgue measure, then $\hat{w}_n(\gamma_n(U_1))$ is uniform on the vertices of $\mathcal{T}_n$. Alternatively, we have that the measure $\mu_n$ defined by, for $A\subseteq\mathcal{T}_n$,
\[\mu_n(A):=\lambda \{t\in[0,1]:\:\hat{w}_n(\gamma_n(t))\in A\},\]
is uniform on the vertices of $\mathcal{T}_n$.

We can now define the sub-trees $\mathcal{T}_n(k)$. For $n\geq 1$, define a collection of vertices by $v_k^n:=\hat{w}_n(\gamma_n(u^n_k))$, and, for $k\geq 1$, let
\[\mathcal{T}_n(k):=r(\mathcal{T}_n,\{v^n_1,\dots,v^n_k\}).\]
The corresponding measure projection of $\mu_n$ onto $\mathcal{T}_n(k)$ is denoted
\begin{equation}\label{munk}
\mu_n^{(k)}:=\mu_n\circ\phi_{\mathcal{T}_n,\mathcal{T}_n(k)}^{-1},
\end{equation}
where the projection operator $\phi_{\mathcal{T}_n,\mathcal{T}_n(k)}$ is defined on graph trees analogously to the projection operator for dendrites, (see (\ref{project})).

One of the main results in \cite{Aldous3} is Theorem 20, in which necessary and sufficient conditions for the distributional convergence of rescaled search-depth functions to the normalised Brownian excursion are presented. These are the convergence of random finite dimensional distributions and a tightness result. We now translate one part of this result into our setting, although we omit the proof since it may be demonstrated by repeating exactly the same steps as were used in the proof of \cite{Aldous3}, Theorem 20. Analogous to the definition of $\Delta^{(k)}$ in Lemma \ref{deltakdecay}, we introduce the notation
\begin{equation}\label{deltankdef}
\Delta_n^{(k)}:=\sup_{x\in\mathcal{T}_n}d_{\mathcal{T}_n}\left(x,\phi_{\mathcal{T}_n,\mathcal{T}_n(k)}(x)\right).
\end{equation}

{\lem \label{equiv} Under Assumption 1, we have that $\lim_{k\rightarrow\infty}\limsup_{n\rightarrow\infty}n^{-1/2}\Delta_n^{(k)}=0$.}
\bigskip

Finally, we note that using the sequential construction of \cite{Aldous3}, Section 2.2, which was outlined briefly in Section \ref{crtsec}, for each $n$, we can isometrically embed the vertices of $\mathcal{T}_n$ into $l^1$ from the vertex sequence $(v^n_k)_{k\geq 1}$. Observe that, under Assumption 1, because we are assuming $(u^n_k)_{k\geq 1}$ to be dense in $[0,1]$, the sequence $(v^n_k)_{k\geq 1}$ will contain all the vertices of $\mathcal{T}_n$, and so this procedure does result in an isometric embedding for $(\mathcal{T}_n,\mu_n)$. We shall denote by $\psi_n$ the distance-preserving map from the vertices of $\mathcal{T}_n$ into $l^1$, and by $\tilde{\mathcal{T}_n}, \tilde{\mu}_n,\dots$ the $l^1$ embedded versions of objects. As with the embeddings for dendrites, we will discuss the measurability of this procedure in Section \ref{meassec}.

\subsection{Discrete processes}\label{discreteproc}

We now define the various discrete processes that appear in this article. Here and elsewhere we apply the convention that the notation $m$ represents a discrete time parameter, in contrast to the continuous time parameter $t$. Throughout this section, we shall assume that we have been given a fixed realisation of $\mathcal{T}_n$ and $(\mathcal{T}_n(k))_{k\geq 1}$.

The fundamental process of interest is $(X^n_m)_{m\geq 0}$, the simple random walk on the vertices of $\mathcal{T}_n$, which we shall suppose is built on our underlying probability space. By simple random walk, we mean the process describing the position of a particle, started from the root, which jumps at each time step to a neighbouring vertex of $\mathcal{T}_n$, with equal probability being placed on each of the possible choices, and each jump being independent of the past (apart from the position of the particle at that time). The law of $X^n$ will be denoted $\mathbf{P}^{\mathcal{T}_n}_\rho$, and its image in $l^1$ under the distance-preserving map $\psi_n$ introduced at the end of the previous section (when this is defined) by $\tilde{\mathbf{P}}^{\mathcal{T}_n}_\rho:=\mathbf{P}^{\mathcal{T}_n}_\rho\circ\psi_n^{-1}$. As in the introduction, we extend $\tilde{\mathbf{P}}^{\mathcal{T}_n}_\rho$ to a law on the continuous paths in $l^1$ by linear interpolation of discrete sample paths.

The first related process we construct is simply the projection of the simple random walk $X^n$ onto $\mathcal{T}_n(k)$, we shall denote this by $(X^{n,k}_m)_{m\geq 0}$, and define it precisely by
\begin{equation}\label{xnkdef}
X^{n,k}_m:=\phi_{\mathcal{T}_n,\mathcal{T}_n(k)}(X^n_m).
\end{equation}
The associated jump process we shall write as $(J^{n,k}_m)_{m\geq 0}$. Of course, $J^{n,k}$ is nothing more than the simple random walk on the vertices of $\mathcal{T}_n(k)$. It will be useful to be able to express $X^{n,k}$ in terms of $J^{n,k}$, and to do this we introduce a process $({A}^{n,k}_m)_{m\geq 0}$ that is defined by $A^{n,k}_0=0$ and, for $m\geq 1$,
\begin{equation}\label{ankdef}
A_m^{n,k}:=\min\{l\geq A_{m-1}^{n,k}:\:X^n_{l}\in\mathcal{T}_n(k)\backslash\{X^{n}_{A_{m-1}^{n,k}}\}\},
\end{equation}
so that the time $A_m^{n,k}-A_{m-1}^{n,k}$ is the time until the random walk $X^n$ hits a vertex in $\mathcal{T}_n(k)$ other than the one it was in at time $A_{m-1}^{n,k}$. If we then define $(\tau^{n,k}(m))_{m\geq 0}$ by
\begin{equation}\label{taudef}
{\tau}^{n,k}(m):=\max\{l:\:{A}^{n,k}_{l}\leq m\},
\end{equation}
it is easy to check that $X^{n,k}$ is recovered by taking
\begin{equation}\label{recover}
X^{n,k}_m=J^{n,k}_{\tau^{n,k}(m)}.
\end{equation}

One of the key steps in our proof of Theorem \ref{mainresult} is showing that, as $n$ and then $k$ become large, the process $A^{n,k}$ may be rescaled to a function that is linear in time, see Corollary \ref{agrowth}. However, the process $A^{n,k}$ is relatively difficult to handle directly, and so we now introduce a closely related process that is more manageable. First, we define the occupation times, $(\ell^{n,k}_m(x))_{m\geq 0,\:x\in\mathcal{T}_n(k)}$, and a stationary measure, $\nu_n^{(k)}$, of the jump process $J^{n,k}$ by setting
\begin{equation}\label{occtime}
\ell^{n,k}_m(x):=\sum_{l=0}^{m}\mathbf{1}_x(J_{l}^{n,k}),\hspace{20pt}\nu_n^{(k)}(\{x\}):=\frac{\mathrm{deg}_{n,k}(x)}{2},
\end{equation}
for  $x$ a vertex in $\mathcal{T}_n(k)$, where $\mathrm{deg}_{n,k}:=\mathrm{deg}_{\mathcal{T}_n(k)}$. From these quantities we define the local times (or occupation time densities), $(L^{n,k}_m(x))_{m\geq 0,\:x\in\mathcal{T}_n(k)}$, of the jump process by
\begin{equation}\label{loctime}
L_m^{n,k}(x):=\frac{\ell^{n,k}_m(x)}{\nu_n^{(k)}(\{x\})}.
\end{equation}
We can use these local times to define an additive functional, $(\hat{A}^{n,k}_m)_{m\geq 0}$, by $\hat{A}^{n,k}_0=0$, and for $m\geq1$,
\begin{equation}\label{tildeankdef}
\hat{A}^{n,k}_m:=n\int_{\mathcal{T}_n(k)}L^{n,k}_{m-1}(x)\mu_n^{(k)}(dx).
\end{equation}
As a result of this integral representation, $\hat{A}^{n,k}$ is much easier to deduce convergence results for than $A^{n,k}$; it also gives a good approximation of $A^{n,k}$. The reason for this second fact is explained by the following. Extend $\hat{A}^{n,k}$ to continuous time by linear interpolation and let the (continuous time) inverse of $\hat{A}^{n,k}$ be defined by $\hat{\tau}^{n,k}(t):=\max\{s:\hat{A}^{n,k}_s\leq t\}$. Now introduce a time-changed version of $J^{n,k}$, denoted $(\hat{X}^{n,k}_t)_{t\geq 0}$, and defined by
\begin{equation}\label{hatx}
\hat{X}^{n,k}_t:=J^{n,k}_{\hat{\tau}^{n,k}(t)}.
\end{equation}
Clearly, both $X^{n,k}$ and $\hat{X}^{n,k}$ have by construction the same jump chain, $J^{n,k}$. The process $X^{n,k}$ sits at a vertex in $\mathcal{T}_n(k)$ while $X^n$ jumps about in $\mathcal{T}_n$ until $X^n$ hits a different vertex in $\mathcal{T}_n(k)$, and so the length of time spent in each place is a (possibly unbounded) random variable. The process $\hat{X}^{n,k}$, on the other hand, waits at a vertex $x$ a fixed time $2n\mu_n^{(k)}(\{x\})/\mathrm{deg}_{n,k}(x)$ before jumping. The processes $X^{n,k}$ and $\hat{X}^{n,k}$ can be shown to be close when suitably rescaled, and the reason for this is that the time $2n\mu_n^{(k)}(\{x\})/\mathrm{deg}_{n,k}(x)$ gives a good approximation of the expectation of the random time that $X^{n,k}$ must wait at vertex $x$ before jumping. More specifically, we prove a tightness result for $A^{n,k}$ and $\hat{A}^{n,k}$, see Proposition \ref{anktildeank}.

\subsection{Overview of proof}\label{overv}

As with any long proof, there is a danger that the main arguments will be lost in amongst the details and technicalities. To try to avoid this problem, we present here a brief summary of the key steps, and an index of processes is provided in Appendix \ref{notation}. Pictorially, we have that the processes are related in the following fashion
\[n^{-1/2}X^n_{tn^{3/2}} \hspace{8pt} \approx\hspace{8pt} n^{-1/2}X^{n,k}_{tn^{3/2}}\hspace{8pt} \approx\hspace{8pt} n^{-1/2}\hat{X}^{n,k}_{tn^{3/2}}=n^{-1/2}J^{n,k}_{\hat{\tau}^{n,k}(tn^{3/2})}\hspace{8pt} \rightsquigarrow\hspace{8pt} B^{(k)}_t\hspace{8pt} \rightarrow\hspace{8pt} X_t.\]

The process $X^{n,k}$ is the projection of $X^n$ onto $\mathcal{T}_n(k)$, and so to prove that the two processes are close, we need to show that the projection operator $\phi_{\mathcal{T}_n,\mathcal{T}_n(k)}$ does not move points too far. This purely geometrical result, which is stated as Lemma \ref{lemyo}, is covered by Lemma \ref{equiv}. The connection between $X^{n,k}$ and $\hat{X}^{n,k}$ was discussed at the end of the previous section.

The point of transfer between discrete and continuous time processes is Proposition \ref{jnkbkl1}, where we demonstrate the unsurprising result that, when rescaled, the simple random walks $J^{n,k}$ on $\mathcal{T}_n(k)$ converge as $n\rightarrow\infty$ to the Brownian motion $B^{(k)}$ on $\mathcal{T}{(k)}$. In showing that the limit of $\hat{X}^{n,k}$ is also close to $B^{(k)}$, by applying the representation at (\ref{hatx}) it will suffice to exhibit the behaviour of $\hat{A}^{n,k}$ as $n$ and then $k$ gets large. The two concrete results we prove are the following. Firstly, by demonstrating that the rescaled local times of the jump processes $J^{n,k}$ converge when rescaled to those of $B^{(k)}$, (see Lemma \ref{localtimeconv}), we are able show that $\hat{A}^{n,k}$, as defined at (\ref{tildeankdef}), may be rescaled to converge to a related additive functional, $\hat{A}^{(k)}$, defined from $B^{(k)}$, (see (\ref{tildeakdef}) for a definition of $\hat{A}^{(k)}$ and Corollary \ref{tildeankak} for a statement of the relevant convergence result). Secondly, we deduce that  $\hat{A}^{(k)}_t$ converges to $t$ uniformly on compact intervals (Proposition \ref{tildeakt}). Although we will not proceed to present these results rigourously in the way we now describe, the motivation for our argument is provided by the following. First, it is possible to deduce that
\[n^{-1/2}\hat{X}^{n,k}_{tn^{3/2}}=n^{-1/2}J^{n,k}_{\hat{\tau}^{n,k}(tn^{3/2})}\rightarrow B^{(k)}_{\hat{\tau}^{(k)}(t)},\]
where $\hat{\tau}^{(k)}(t)$ is the right continuous inverse of $\hat{A}^{(k)}$, defined similarly to (\ref{taukdef}). Since $\hat{A}^{(k)}_t\rightarrow t$ uniformly on compact intervals, $\hat{\tau}^{(k)}(t)\rightarrow t$ on compact intervals. Thus the continuity of $B^{(k)}$ implies that $B^{(k)}_{\hat{\tau}^{(k)}(t)}$ is close to $B^{(k)}_t$ uniformly.

The construction of $B^{(k)}$ in Section \ref{crtproc} using a time-change argument allows us to prove an almost-sure version of the limit $B^{(k)}\rightarrow X$ via standard arguments, which depend only on the fact that the measures $\lambda^{(k)}$ converge weakly to $\mu$, see Lemma \ref{bkx}.

\section{Convergence of Brownian motion on finite trees}

In this section, we fix $(w,u)\in\Gamma$, and show that if the processes $(B^{(k)})_{k\geq 1}$ and $X$ are coupled as in Section \ref{crtproc}, then $B^{(k)}$ converges $\mathbf{P}$-a.s. on any compact time interval to $X$ as $k\rightarrow\infty$. We will also prove the convergence of a related additive functional.

{\lem \label{bkx} Fix $(w,u)\in\Gamma$ and $R\in(0,\infty)$. If the processes $(B^{(k)})_{k\geq 1}$ and $X$ are coupled as in Section \ref{crtproc}, then $\mathbf{P}$-a.s., \[(B^{(k)}_t)_{t\in[0,R]}\rightarrow(X_t)_{t\in[0,R]},\hspace{20pt}\mbox{as }k\rightarrow\infty,\]
in $C([0,R],\mathcal{T})$.}
\begin{proof} We start by demonstrating that, $\mathbf{P}$-a.s.,
\begin{equation}\label{unicon}
\sup_{t\in[0,R+1]}|{A}^{(k)}_t-t|\rightarrow0,\hspace{20pt}\mbox{as }k\rightarrow\infty,
\end{equation}
where ${A}^{(k)}$ is the additive functional defined at (\ref{akdef}). Fix $\varepsilon>0$. By Lemma \ref{localtimecont} and the definition of $\Gamma$, we can assume that the local times of $X$ are jointly continuous and $\lambda^{(k)}$ converges weakly to $\mu$. These assertions imply that, point-wise for $t\geq0$, we have
\[{A}_t^{(k)}=\int_{\mathcal{T}(k)}L_{t}(x)\lambda^{(k)}(dx)\rightarrow\int_{\mathcal{T}}L_{t}(x)\mu(dx)=t.\]
Note that the integral over $\mathcal{T}(k)$ makes sense because, by construction, $\mathcal{T}(k)\subseteq\mathcal{T}$. Using the monotonicity in $t$ of the functions $A^{(k)}$, we can apply an elementary argument to deduce from this that the uniform convergence at (\ref{unicon}) holds.

As a consequence of (\ref{unicon}), we also have that $\sup_{t\in[0,R]}|\tau^{(k)}_t-t|\rightarrow0$, where $\tau^{(k)}$ is the inverse of $A^{(k)}$ defined at (\ref{taukdef}). Recalling from (\ref{bkdef}) that $B^{(k)}_t:=X_{\tau^{(k)}(t)}$, the $\mathbf{P}$-a.s. continuity of $X$ implies the result.
\end{proof}

To state the corresponding result for convergence of probability laws in $l^1$, we introduce the notation
\begin{equation}\label{tpk}
\tilde{\mathbf{P}}^{\mathcal{T}(k)}_\rho:=\mathbf{P}^{\mathcal{T}(k),\lambda^{(k)}}_\rho\circ\psi^{-1},
\end{equation}
where $\psi:\mathcal{T}\rightarrow l^1$ is the distance-preserving map introduced at the end of Section \ref{crtsec}.

{\propn \label{prop1} If $(w,u)\in\Gamma$, then
\[\left(\tilde{\mathcal{T}}(k),\tilde{\mu}^{(k)},\tilde{\mathbf{P}}^{\mathcal{T}(k)}_\rho\right)\rightarrow\left(\tilde{\mathcal{T}},\tilde{\mu},\tilde{\mathbf{P}}^{\mathcal{T}}_\rho\right),\]
in the space $\mathcal{K}(l^1)\times \mathcal{M}_1(l^1)\times\mathcal{M}_1(C([0,1],l^1))$.}
\begin{proof} It is easy to check from the construction of $\tilde{\mathcal{T}}(k)$ and $\tilde{\mu}^{(k)}$ that both $d_H^{l^1}(\tilde{\mathcal{T}}(k),\tilde{\mathcal{T}})$ and $d_P^{l^1}(\tilde{\mu}^{(k)},\tilde{\mu})$ are bounded above by $\Delta^{(k)}$, where $d_H^{l^1}$ is the Hausdorff metric on $\mathcal{K}(l^1)$, $d_P^{l^1}$ is the Prohorov metric on $\mathcal{M}_1(l^1)$, and $\Delta^{(k)}$ is the quantity defined in Lemma \ref{deltakdecay}. Applying Lemma \ref{deltakdecay} we immediately are able to deduce that $(\tilde{\mathcal{T}}(k),\tilde{\mu}^{(k)})\rightarrow(\tilde{\mathcal{T}},\tilde{\mu})$ in the appropriate space.

Define now $\tilde{B}^{(k)}:=\psi(B^{(k)})$ and $\tilde{X}:=\psi(X)$, where $B^{(k)}$ and $X$ are coupled as in Section \ref{crtproc}. Applying the fact that $\psi$ is distance-preserving and Lemma \ref{bkx} , we have that $\tilde{B}^{(k)}\rightarrow \tilde{X}$, $\mathbf{P}$-a.s., in $C([0,1],l^1)$. Since $\tilde{B}^{(k)}$ has law $\tilde{\mathbf{P}}^{\mathcal{T}(k)}_\rho$ and $\tilde{X}$ has law $\tilde{\mathbf{P}}^{\mathcal{T}}_\rho$, the result follows.
\end{proof}

That local times of $B^{(k)}$ exist is guaranteed by the following lemma.

{\lem \label{localtimekcont} Fix $(w,u)\in\Gamma$. If $B^{(k)}$ is a random process with law $\mathbf{P}^{\mathcal{T}(k),\lambda^{(k)}}_\rho$, then $B^{(k)}$ admits local times $(L^{(k)}_t(x))_{t\geq 0,\:x\in\mathcal{T}(k)}$ that are $\mathbf{P}$-a.s. jointly continuous in $t$ and $x$.}
\begin{proof} The existence and continuity of local times for $B^{(k)}$ may be shown in exactly the same way as for the process $X$, see Lemma \ref{localtimecont}.  However, to do this, it is necessary to provide suitable estimates for the size of an $\varepsilon$-cover for $\mathcal{T}(k)$ and on the heat kernel of $B^{(k)}$ in place of (\ref{coversize}) and (\ref{tdupper}), respectively. First, since $\mathcal{T}(k)$ is made up of a finite collection of line segments, and $\lambda^{(k)}$ is simply the rescaled Lebesgue measure on these, there is no difficulty in deducing that there exist constants $c_1$, $c_2$ and $r_0>0$ such that
\[c_1r \leq \lambda^{(k)}(B_{\mathcal{T}(k)}(x,r))\leq c_2 r,\hspace{20pt}\forall x\in\mathcal{T}(k),\:r\in(0,r_0),\]
where $B_{\mathcal{T}(k)}(x,r)$ is the ball of radius $r$ around $x$ in $\mathcal{T}(k)$. This allows us to apply \cite{Kumagai}, Theorem 3.1, to deduce the existence  of a heat kernel $(p_t^{(k)}(x,y))_{t\geq 0,\:x,y\in\mathcal{T}(k)}$ for $B^{(k)}$ which satisfies, for some $c_3$ and $t_0>0$, $p^{(k)}_t(x,x)\leq c_3 t^{-1/2}$, for all $x\in\mathcal{T}(k)$, $t\in(0,t_0)$. Secondly, we can use again the simple structure of $(\mathcal{T}(k),\lambda^{(k)})$ to deduce that there exists a constant $c_4$ such that $N(\mathcal{T}(k),\varepsilon)\leq c_4 \varepsilon^{-1}$, for every $\varepsilon\in(0,1)$. These two estimates enable us to complete the proof using the argument of Lemma \ref{localtimecont}.
\end{proof}

We now introduce another additive functional, $(\hat{A}^{(k)}_t)_{t\geq0}$, that we will later show describes the scaling limit as $n\rightarrow\infty$ of the functions $\hat{A}^{n,k}$, as defined at (\ref{tildeankdef}). Set
\begin{equation}\label{tildeakdef}
\hat{A}^{(k)}_t:=\int_{\mathcal{T}(k)}L^{(k)}_t(x)\mu^{(k)}(dx).
\end{equation}
The following description of the local times of $B^{(k)}$ will be useful in demonstrating that the additive functionals $\hat{A}^{(k)}$ converge in the subsequent lemma.

{\lem \label{loctimerep} Fix $(w,u)\in\Gamma$. If the processes $(B^{(k)})_{k\geq 1}$ and $X$ are coupled as in Section \ref{crtproc}, then $\mathbf{P}$-a.s., the local times $(L^{(k)}_t(x))_{t\geq 0,\:x\in\mathcal{T}(k)}$ of $B^{(k)}$ satisfy
\[L_t^{(k)}(x)=L_{{\tau}^{(k)}(t)}(x),\hspace{20pt}\forall t\geq0,\:x\in\mathcal{T}(k),\]
where $(L_t(x))_{t\geq 0, x\in\mathcal{T}}$ are the local times of $X$.}
\begin{proof} The following argument holds $\mathbf{P}$-a.s. Fix $k\geq 1$. Assuming that ${L}^{(k)}_t(x)$ is jointly continuous in $t$ and $x$, it is possible to deduce that the maps
\[A\mapsto\int_{A}d{A}_u^{(k)},\hspace{20pt}A\mapsto\int_{\mathcal{T}(k)} \int_{A} dL_u(x) \lambda^{(k)}(dx),\]
for Borel sets $A\subseteq \mathbb{R}_+$, are well-defined and describe Borel measures on $\mathbb{R}_+$. Furthermore, for an interval $(s,t]\subseteq \mathbb{R}_+$, we have
\[\int_{(s,t]}d{A}_u^{(k)}={A}_t^{(k)}-{A}_s^{(k)}=\int_{\mathcal{T}(k)} (L_t(x)-L_s(x))\lambda^{(k)}(dx)=\int_{\mathcal{T}(k)} \int_{(s,t]} dL_u(x) \lambda^{(k)}(dx).\]
By a standard argument (see \cite{Kallenberg}, Theorem 2.14, for example) if two locally finite Borel measures on $\mathbb{R}_+$ agree on sets of the form $(s,t]$ and have no atom at zero, they are identical. Applying this fact, for a measurable $B\subseteq\mathcal{T}(k)$, we have
\[\int_{0}^t\mathbf{1}_B(B_s^{(k)})ds=\int_0^{{\tau}^{(k)}(t)}\mathbf{1}_B(X_s)d{A}_s^{(k)}=\int_{\mathcal{T}(k)}\int_{[0,{\tau}^{(k)}(t)]\cap X^{-1}(B)} dL_u(x) \lambda^{(k)}(dx),\]
where
$X^{-1}(B):=\{s:\:X_s\in B\}$ is a measurable subset of $\mathbb{R}_+$. An elementary argument using the continuity of $X$ allows it to be deduced that the measure $dL_u(x)$ is supported on the set $X^{-1}(\{x\})$. Hence
\[\int_{0}^t\mathbf{1}_B(B_s^{(k)})ds=\int_{\mathcal{T}(k)}\int_{[0,{\tau}^{(k)}(t)]} \mathbf{1}_B(x) dL_u(x) \lambda^{(k)}(dx)=\int_{\mathcal{T}(k)}\mathbf{1}_B(x) L_{{\tau}^{(k)}(t)}(x) \lambda^{(k)}(dx),\]
from which the result follows.\end{proof}

{\propn \label{tildeakt} Fix $(w,u)\in\Gamma$ and $R>0$. If the processes $(B^{(k)})_{k\geq 1}$ and $X$ are coupled as in Section \ref{crtproc}, then $\mathbf{P}$-a.s.,
\[\sup_{t\in[0,R]}\left|\hat{A}^{(k)}_t-t\right|\rightarrow 0,\hspace{20pt}\mbox{as }k\rightarrow \infty.\]}
\begin{proof} The following proof holds $\mathbf{P}$-a.s. By the previous lemma and the definition of $\hat{A}^{(k)}$, we have that, for $t\geq 0$,
\[\hat{A}^{(k)}_t=\int_{\mathcal{T}(k)}L_{{\tau}^{(k)}(t)}(x)\mu^{(k)}(dx)= \int_{\mathcal{T}}L_{{\tau}^{(k)}(t)}(\phi_{\mathcal{T},\mathcal{T}(k)}(x))\mu(dx),\]
where for the second equality we use the definition of $\mu^{(k)}$ as the projection of $\mu$ onto $\mathcal{T}(k)$. It immediately follows that
\[\sup_{t\in[0,R]}\left|\hat{A}^{(k)}_t - t\right|\leq \sup_{t\in[0,R]}\sup_{\buildrel{\scriptstyle {x,y\in\mathcal{T}:}}\over {d_{\mathcal{T}}(x,y)\leq \Delta^{(k)}}} \left| L_{{\tau}^{(k)}(t)}(x)-L_{{\tau}^{(k)}(t)}(y)\right|+\sup_{t\in[0,R]}\left|{\tau}^{(k)}(t)-t\right|.\]
In the proof of Lemma \ref{bkx}, we showed that ${\tau}^{(k)}(t)\rightarrow t$ uniformly on $[0,R]$. Combining this result with the fact that $\Delta^{(k)}\rightarrow 0$ and the continuity of the local times of $X$ (see Lemmas \ref{deltakdecay} and \ref{localtimecont}, respectively), it is straightforward to use the above estimate derive the result.
\end{proof}

\section{Convergence of jump processes and local times}\label{finiteconverge}

The primary aim of this section is to demonstrate that the processes $J^{n,k}$, when rescaled, converge in distribution to $B^{(k)}$ as $n\rightarrow\infty$. We also show that the additive functional $\hat{A}^{n,k}$ defined at (\ref{tildeankdef}) converges to the process $\hat{A}^{(k)}$ introduced in the previous section.  A key result is Proposition \ref{corr}, where we show the simultaneous convergence of trees, measures, jump processes and local times, and from which the convergence of $\hat{A}^{n,k}$ follows easily using the continuous mapping theorem. For the purposes of this section, because the trees we discuss have a finite number of branches, it will be convenient to work in the space of abstract trees with edge lengths using the topology we now introduce.

We consider elements of the form $(T, \mu, f_1, f_2)$. Here, ${T}=({T}^*;|e_1|,\dots,|e_{l-1}|)$ for some $l$, where ${T}^*$ is an ordered graph tree with $l$ vertices, and $|e_1|,\dots,|e_{l-1}|$ are the edge lengths. By including line segments along edges, naturally associated with ${T}$ is a dendrite $\underline{T}$ equipped with the natural shortest path metric $d_{\underline{T}}$. We assume that $\mu$ is a Borel probability measure on $\underline{T}$. The continuous $\underline{T}$-valued function $f_1$ is defined on some interval $[0,R]$, and can be considered as the sample path of a process on $\underline{T}$. Finally, the $\mathbb{R}_+$-valued function $f_2$ is defined on $[0,R]\times \underline{T}$ and can be thought of as representing the corresponding local times.

To define the topology of interest, we introduce a metric, $d$, between two such 4-tuplets, $(T, \mu, f_1,f_2)$ and $(T',\mu', f_1',f_2')$, when the intervals on which the functions $f_1$, $f_1'$ are defined are the same. First, we introduce a distance $d_1$ between ordered graph trees with edge lengths. If ${T}^*\neq {T}'^*$, then set $d_1(T,T')=\infty$. Otherwise, assume ${T}^*={T}'^*$. The distance between trees is defined to be the maximal edge length difference, i.e.,
\[d_1(T,T'):=\sup_{i}\left| |e_i|-|e_i'|\right|.\]
When ${T}^*={T}'^*$, we have a homeomorphism $\Upsilon_{\underline{T},\underline{T'}}:\underline{T}\rightarrow \underline{T'}$, under which the point $x\in \underline{T}$, which is a distance $\alpha$ along the edge $e_i$ (considered from the vertex at the end of $e_i$ which is closest to the root), is mapped to the point $x'\in \underline{T'}$ which is a distance $|e_i'|\alpha/|e_i|$ along $e_i'$. We use this function to define a collection of distances. Let
\[d_2(\mu,\mu'):=d_P(\underline{T};\mu,\mu'\circ\Upsilon_{\underline{T},\underline{T'}})+d_P(\underline{T'};\mu\circ\Upsilon_{\underline{T'},\underline{T}},\mu'),\]
where $d_P(\underline{T};\cdot,\cdot)$ is the usual Prohorov metric on $\underline{T}$, and we make this choice for the reason that it induces the weak topology on $\underline{T}$. Furthermore, set
\[d_3(f_1,f_1'):=\sup_{t\in[0,R]} \left[ d_{\underline{T}} (f_1(t),\Upsilon_{\underline{T'},\underline{T}}(f_1'(t)))+ d_{\underline{T'}} (\Upsilon_{\underline{T},\underline{T'}}(f_1(t)),f_1'(t)) \right],\]
\[d_4(f_2,f_2'):=\sup_{t\in[0,R],\:x\in \underline{T}}\left|f_2(t,x)-f_2'(t,\Upsilon_{\underline{T},\underline{T'}}(x))\right|.\]
The metric $d$ is then defined by setting,
\[d((T, \mu, f_1,f_2),(T', \mu' f_1',f_2')):=\left(d_1(T,T')+d_2(\mu,\mu')+d_3(f_1,f_1')+d_4(f_2,f_2')\right)\wedge 1.\]

We are now almost in a position to state and prove the first result of this section. In this lemma and subsequent results of this section, we assume that we have been given the collections  $(\mathcal{T}_n)_{n\geq 1}$ and $(\mathcal{T}_n(k))_{n,k\geq1}$, and that these are constructed from a sequence $\{(w_n,u^n)\}_{n\geq 1}$ that satisfies Assumption 1. We define $T_n(k)$ to be the graph tree with vertices given by the root and leaves of the graph tree $\mathcal{T}_n(k)$ along with their branch points in $\mathcal{T}_n(k)$. The edge lengths of $T_n(k)$ are those induced from the graph distance $d_{\mathcal{T}_n}$ on $\mathcal{T}_n(k)$, and the ordering of vertices of ${T}_n(k)$ follows from the ordering of vertices of $\mathcal{T}_n$. Since there is a natural distance-preserving embedding of the vertices of the graph tree $\mathcal{T}_n(k)$ into the dendrite $\underline{T}_n(k)$, the measure $\mu_n^{(k)}$, as described at (\ref{munk}), may be thought of as a Borel measure on $\underline{T}_n(k)$ consisting of a finite number of atoms. Similarly, we define $T(k)$ to be the ordered graph tree with edge lengths constructed from the pair $(w,u)\in\Gamma$ that corresponds to the limit of $\{(w_n,u^n)\}_{n\geq 1}$. Also, the measure $\mu^{(k)}$ can be thought of as a Borel measure on $\underline{T}(k)$, which is a dendrite with exactly the same structure as $\mathcal{T}(k)$. Finally, we also introduce notation for the rescaled trees and measures, specifically, we set
\[\breve{T}_n(k):=({T}_n(k)^*; n^{-1/2}|e_1|,\dots,n^{-1/2}|e_{l-1}|),\]
where $T_n(k)$ has $l$ vertices. We define $\breve{\mu}_n^{(k)}$ to be the probability measure on the dendrite associated with $\breve{{T}}_n(k)$ satisfying
\[\breve{\mu}_n^{(k)}=\mu_n^{(k)}\circ \Upsilon_{\breve{\underline{T}}_n(k),\underline{T}_n(k)}.\]

{\lem \label{subtreelemma} Under Assumption 1, $(\breve{{T}}_n(k),\breve{\mu}_n^{(k)})\rightarrow ({T}(k),\mu^{(k)})$, as $n\rightarrow\infty$, with respect to the distance $d_1+d_2$.}
\begin{proof} The result that $d_1(\breve{{T}}_n(k),{T}(k))\rightarrow 0$ is essentially demonstrated in the proof of \cite{Aldous3}, Theorem 20, and so we will restrict ourselves to showing that $d_2(\breve{\mu}_n^{(k)},\mu^{(k)})\rightarrow0$. First, denote by $\varphi_k$ the map from $[0,1]$ to $\underline{T}(k)$ that is obtained by composing the map $t\mapsto [t]$, where $[t]$ is the equivalence class of $t\in[0,1]$, as defined by (\ref{equivrel}), with the projection map $\phi_{\mathcal{T},\mathcal{T}(k)}$ (here, we identify $\mathcal{T}(k)$ and $\underline{T}(k)$ in the obvious way). The (non-root) leaves of $T(k)$ are described by the points $\zeta_i:=\varphi_k(u_i)$, $i=1,\dots,k$, and we shall denote by $b_{ij}$ the branch point of $\rho$, $\zeta_i$ and $\zeta_j$ in $\underline{T}(k)$.

For $t\in[0,1]$, we must have that $t\in[u_i,u_j]$ for some $i,j\in\{-1,0,\dots,k\}$, where we introduce the notation $u_{-1 }=0$ and $u_0=1$, and we assume that $(u_i,u_j)\cap\{u_{-1},u_0,\dots,u_k\}=\emptyset$. A simple analysis of the construction of $\underline{T}_n(k)$ allows it to be deduced that
\begin{equation}\label{defn}
\varphi_{k}(t)=\left\{\begin{array}{ll}
                                                             \mbox{$[[$}b_{ij},\zeta_i\mbox{$]]$}(m_w(t,t_i)-m_w(t_i,t_j)),& \mbox{if }m_w(t,t_i)\geq m_w(t,t_j), \\
                                                               \mbox{$[[$}b_{ij},\zeta_j\mbox{$]]$}(m_w(t,t_j)-m_w(t_i,t_j)),  & \mbox{otherwise,}
                                                              \end{array}
\right.\end{equation}
where $m_w$ is the minimum function defined in Section \ref{abstracttree}, and we use the notation $[[b_{ij},\zeta_i]]\alpha$ to represent the point of $\underline{T}_n(k)$ that lies on the line segment $[[b_{ij},\zeta_i]]$ a distance $\alpha$ from $b_{ij}$. Note that the right hand-side of the above expression is well-defined if we set $\zeta_0=\zeta_{-1}=\rho$. Furthermore, observe that if $m_w(t,t_i)= m_w(t,t_j)$ then the two expressions in the right-hand side of (\ref{defn}) are equal (to $b_{ij}$).

Analogous to the above definition, we set $\varphi_{n,k}(t):=\Upsilon_{\underline{{T}}_n(k),\breve{\underline{{T}}}_n(k)}(\phi_{\mathcal{T}_n,\mathcal{T}_n(k)}(\hat{w}_n(\gamma_n(t))))$, which is a map from $[0,1]$ to $\breve{\underline{T}}_n(k)$ (we consider that vertices of $\mathcal{T}_n(k)$ are embedded in $\underline{T}_n(k)$ in the natural way). Also, denote $\zeta^n_i:=\varphi_{n,k}(u^n_i)$, and the branch point of $\rho$, $\zeta^n_i$ and $\zeta^n_j$ in $\breve{\underline{T}}_n(k)$ by $b^n_{ij}$. An expression for $\varphi_{n,k}$ of the form of (\ref{defn}) is not difficult to deduce.

Since $d_1(\breve{{T}}_n(k),{T}(k))\rightarrow 0$, for large $n$ we can define the homeomorphism $\Upsilon_{\breve{\underline{T}}_n(k),\underline{T}(k)}$ from $\breve{\underline{T}}_n(k)$ to $\underline{T}(k)$ by rescaling edges in the way described at the start of this section. For large $n$, we clearly have that $\Upsilon_{\breve{\underline{T}}_n(k),\underline{T}(k)}(\zeta^n_i)=\zeta_i$ and $\Upsilon_{\breve{\underline{T}}_n(k),\underline{T}(k)}(b^n_{ij})=b_{ij}$ for each $i,j\leq k$. Consequently, by considering the expression at (\ref{defn}) and similar formulae for $\varphi_{n,k}$, we have that under Assumption 1, for every $t\in[0,1]$,
\begin{equation}\label{convv}
\Upsilon_{\breve{\underline{T}}_n(k),\underline{T}(k)}(\varphi_{n,k}(t))\rightarrow\varphi_{k}(t).
\end{equation}
By definition, $\breve{{\mu}}_n^{(k)}=\lambda\circ\varphi_{n,k}^{-1}$ and $\mu^{(k)}=\lambda\circ\varphi_k^{-1}$, where $\lambda$ is the usual Lebesgue measure on $[0,1]$. Hence, applying the convergence at (\ref{convv}) and Fatou's lemma, for open $A\subseteq\underline{T}(k)$ it follows that $\liminf_{n\rightarrow\infty}\breve{\mu}_n^{(k)}\circ\Upsilon_{\underline{{T}}(k),\breve{\underline{T}}_n(k)}(A)\geq{\mu}^{(k)}(A)$, which implies that $\breve{\mu}_n^{(k)}\circ\Upsilon_{\underline{{T}}(k),\breve{\underline{T}}_n(k)}$ converges weakly to $\mu^{(k)}$ as measures on $\underline{T}(k)$, (see \cite{Bill2}, Theorem 2.1). In particular, we have that $d_P\left(\underline{T}(k); \breve{\mu}_n^{(k)}\circ\Upsilon_{\underline{{T}}(k),\breve{\underline{T}}_n(k)},\mu^{(k)}\right)\rightarrow 0$. Finally, the map $\Upsilon_{\underline{{T}}(k),\breve{\underline{T}}_n(k)}$ is Lipschitz, and if $(c_n)_{n\geq 1}$ represents the associated Lipschitz constants, then it follows from $d_1(\breve{{T}}_n(k),{T}(k))\rightarrow 0$ that $c_n\rightarrow 1$. Consequently
\[d_P\left(\breve{\underline{T}}_n(k); \breve{\mu}_n^{(k)},\mu^{(k)}\circ \Upsilon_{\breve{\underline{T}}_n(k),\underline{T}(k)}\right)\leq c_n d_P\left(\underline{T}(k); \breve{\mu}_n^{(k)}\circ\Upsilon_{\underline{{T}}(k),\breve{\underline{T}}_n(k)},\mu^{(k)}\right)\rightarrow 0,\]
which completes the proof.
\end{proof}

To define the jump process on $\breve{\underline{T}}_n(k)$ that will be the focus of this section, we need to clarify what we mean by a vertex and so we introduce the set
\[V\left(\breve{\underline{T}}_n(k)\right):=\{x\in\breve{\underline{T}}_n(k):\:d_{\breve{\underline{T}}_n(k)}(\rho,x)=mn^{-1/2},\mbox{ for some }m\in\mathbb{N}\cup\{0\}\}\]
to represent the ``rescaled graph vertices'' contained in $\breve{\underline{T}}_n(k)$. We then define the process $\breve{J}^{n,k}:=(\breve{J}^{n,k}_m)_{m\geq 0}$ to be the simple random walk on $V(\breve{\underline{T}}_n(k))$ started from the root, where we suppose that two elements of $V(\breve{\underline{T}}_n(k))$ are joined by an edge if an only if the line segment between them in $\breve{\underline{T}}_n(k)$ contains no other point in $V(\breve{\underline{T}}_n(k))$. We extend the definition of $\breve{J}^{n,k}$ to all $t\geq 0$ by linear interpolation. To prove convergence of the jump-processes we will need to time-scale $\breve{J}^{n,k}$ according to the length of the graph $\breve{T}_n(k)$. We define
\begin{equation}\label{length}
\Lambda_n^{(k)}:=\sum_{i}|e^{n,k}_i|,
\end{equation}
where $|e^{n,k}_i|$ are the edge lengths of $\breve{T}_n(k)$. Clearly, under Assumption 1, the previous result implies that $\Lambda_n^{(k)}\rightarrow \Lambda^{(k)}$, where $\Lambda^{(k)}$ is the total length of $T(k)$, defined similarly to (\ref{length}). In the following results, we use the notation $B^{(k)}$ to represent the Brownian motion on $(\underline{T}(k), \lambda^{(k)})$ started from the root, where $\lambda^{(k)}$, as defined as at (\ref{leb}), is now thought of as a Borel probability measure on $\underline{T}(k)$. Since $\underline{T}(k)$ and $\mathcal{T}(k)$ are equivalent metric spaces, this is consistent with the definition of $B^{(k)}$ used in earlier sections.

{\lem \label{jnkbk} Fix $R\in(0,\infty)$ and $k\in\mathbb{N}$. Under Assumption 1 it is possible to construct $\breve{J}^{n,k}$ and $B^{(k)}$ under the probability measure $\mathbf{P}$ in such a way that, $\mathbf{P}$-a.s.,
\[\left(\breve{{T}}_n(k),(\breve{J}^{n,k}_{tn\Lambda_{n}^{(k)}})_{t\in[0,R]}\right)\rightarrow\left({{T}}(k),(B^{(k)}_t)_{t\in[0,R]}\right),\]
with respect to $d_1+d_3$.}
\begin{proof} From the previous lemma we have that $\breve{{T}}_n(k)\rightarrow {T}(k)$ with respect to $d_1$. Hence for large $n$ we can define the homeomorphism $\Upsilon_{\breve{\underline{T}}_n(k),\underline{T}(k)}$ from $\breve{\underline{T}}_n(k)$ to $\underline{{T}}(k)$ by rescaling edges in the way described at the start of this section. Now let ${\lambda}_n^{(k)}$ be the scaled Lebesgue measure on $\breve{\underline{T}}_n(k)$, so the mass of a line segment is proportional to its length, and it is normalised so that ${\lambda}_n^{(k)}(\breve{\underline{T}}_n(k))=1$. It is clear that ${\lambda}_n^{(k)}\circ\Upsilon_{\underline{T}(k),\breve{\underline{T}}_n(k)}\rightarrow \lambda^{(k)}$ weakly as probability measures on $\underline{T}(k)$.

Let $B^{(k)}$ be the Brownian motion on $(\underline{T}(k),\lambda^{(k)})$ under the probability measure $\mathbf{P}$. By Lemma \ref{localtimekcont}, we can assume that $B^{(k)}$ has jointly continuous local times, $\mathbf{P}$-a.s.,  which we can use to define the Brownian motion on $(\underline{T}(k),{\lambda}_n^{(k)}\circ\Upsilon_{\underline{T}(k),\breve{\underline{T}}_n(k)})$ by a time-change, similar to that used to define $B^{(k)}$ from $X$ at (\ref{bkdef}). By following the argument of Lemma \ref{bkx} and applying the weak convergence of measures that was noted in the previous paragraph, we are able to deduce that if $\tilde{B}^{n,k}$ is the Brownian motion on $(\underline{T}(k),{\lambda}_n^{(k)}\circ\Upsilon_{\underline{T}(k),\breve{\underline{T}}_n(k)})$ obtained by this time-change, then $\mathbf{P}$-a.s., $\tilde{B}^{n,k}\rightarrow {B}^{(k)}$, in $C([0,R],\underline{T}(k))$. By considering the defining properties of Brownian motion on a dendrite, it is easy to check that under $\mathbf{P}$ the process $B^{n,k}:=\Upsilon_{\underline{T}(k),\breve{\underline{T}}_n(k)}(\tilde{B}^{n,k})$ is  Brownian motion on $(\breve{\underline{T}}_n(k),{\lambda}_n^{(k)})$, and it follows from the previous sentence that, $\mathbf{P}$-a.s.,
\begin{equation}\label{one}
\left(\breve{{T}}_n(k),(B^{n,k}_{t})_{t\in[0,R]}\right)\rightarrow\left({{T}}(k),(B^{(k)}_t)_{t\in[0,R]}\right),
\end{equation}
with respect to $d_1+d_3$.

As a consequence of the hitting probability property of a Brownian motion on a dendrite, if we define $h^{n,k}(0):=0$, and, for $m\geq 1$,
\begin{equation}\label{hitting}
h^{n,k}(m):=\inf\left\{t \geq h^{n,k}(m-1):\:B_t^{n,k}\in V(\breve{\underline{T}}_n(k))\backslash B^{n,k}_{h^{n,k}(m-1)}\right\},
\end{equation}
to be the hitting times of vertices of $V(\breve{\underline{T}}_n(k))$ by $B^{n,k}$, then the discrete time process $(B^{n,k}_{h^{n,k}(m)})_{m\geq 0}$ is a version of the simple random walk on the vertices of $V(\breve{\underline{T}}_n(k))$. Hence we can suppose that $\breve{J}^{n,k}$ is defined by
\begin{equation}\label{couple}
\breve{J}^{n,k}_{m}=B^{n,k}_{h^{n,k}(m)},
\end{equation}
for each $m\in\mathbb{N}$. Furthermore, by considering the excursions of $B^{n,k}$ away from vertices of $V(\breve{\underline{T}}_n(k))$, it is possible to show that $(n\Lambda_n^{(k)}(h^{n,k}(m)-h^{n,k}(m-1)))_{m\geq 1}$
are independent and identically distributed, their common distribution being that of the hitting time of $\{\pm 1\}$ by a standard Brownian motion in $\mathbb{R}^1$, started from zero. Note that the scaling factor of $\Lambda_n^{(k)}$ is a result of the normalisation of ${\lambda}_n^{(k)}$. The hitting time of $\{\pm 1\}$ by a standard Brownian motion in $\mathbb{R}^1$, started from zero, has expected value $1$ and finite fourth moment. As a consequence of these facts, we can apply a standard martingale estimate (\cite{Kallenberg}, Proposition 7.16, for example) to deduce that, for $\varepsilon>0$,
\[\mathbf{P}\left(\sup_{m\leq Rn\Lambda_n^{(k)}}\left|h^{n,k}(m)-\frac{m}{n\Lambda_n^{(k)}}\right|\geq\varepsilon\right)\leq \left(\frac{4}{3\varepsilon}\right)^4\mathbf{E}\left|h^{n,k}(\lfloor Rn\Lambda_n^{(k)}\rfloor)-\frac{\lfloor Rn\Lambda_n^{(k)}\rfloor}{n\Lambda_n^{(k)}}\right|^4\leq  c n^{-2},\]
for some constant $c$ that does not depend on $n$. Note that the second inequality here may be deduced by an elementary argument using the fourth moment condition on the random variables of the form $h^{n,k}(m)-h^{n,k}(m-1)$ (see \cite{Bill}, Theorem 6.1, for example). Thus a Borel-Cantelli argument implies that $\mathbf{P}$-a.s.,
\begin{equation}\label{hittingconv}
\sup_{m\leq Rn\Lambda_n^{(k)}}\left|h^{n,k}(m)-\frac{m}{n\Lambda_n^{(k)}}\right|\rightarrow 0.
\end{equation}
Applying this fact, the convergence result at (\ref{one}), and the coupling of $\breve{J}^{n,k}$ and $B^{n,k}$ from (\ref{couple}), the lemma is readily deduced.
\end{proof}

We are now able to present one of the facts needed in our proof of Theorem \ref{mainresult}. Recall that $J^{n,k}$ is the simple random walk on the vertices of $\mathcal{T}_n(k)$ started from $\rho$. We set $\tilde{J}^{n,k}:=\psi_n(J^{n,k})$, where $\psi_n$ is the distance-preserving map from the vertices of $\mathcal{T}_n$ to $l^1$ introduced at the end of Section \ref{cbptree}. We extend the definition of $\tilde{J}^{n,k}$ by linear interpolation.

{\propn \label{jnkbkl1} Suppose that Assumption 1 holds. If we denote by $\tilde{\mathbf{Q}}_\rho^{\mathcal{T}_n(k)}$ the law of $\left(n^{-1/2}\tilde{J}^{n,k}_{tn\Lambda_n^{(k)}}\right)_{t\in[0,1]}$, then
\[\left(n^{-1/2}\tilde{\mathcal{T}}_n(k),\tilde{\mu}_n^{(k)}(n^{1/2}\cdot),\tilde{\mathbf{Q}}_\rho^{\mathcal{T}_n(k)}\right)\rightarrow \left(\tilde{\mathcal{T}}(k),\tilde{\mu}^{(k)},\tilde{\mathbf{P}}_\rho^{\mathcal{T}(k)}\right)
\]
in the space $\mathcal{K}(l^1)\times \mathcal{M}_1(l^1)\times\mathcal{M}_1(C([0,1],l^1))$, where $\tilde{\mathbf{P}}_\rho^{\mathcal{T}(k)}$ was defined at (\ref{tpk}).}
\begin{proof} Clearly, mapping $(\breve{{T}}_n(k),\breve{\mu}_n^{(k)},(\breve{J}^{n,k}_{tn\Lambda_{n}^{(k)}})_{t\in[0,R]})$ into $l^1$ with the sequential construction (using the vertices $(\zeta_i^n)_{i=1}^k$ defined in the proof of Lemma \ref{subtreelemma}) results in a triple which is identical (in distribution) to $(n^{-1/2}\tilde{\mathcal{T}}_n(k),\tilde{\mu}_n^{(k)}(n^{1/2}\cdot),(n^{-1/2}\tilde{J}^{n,k}_{tn\Lambda_n^{(k)}})_{t\in[0,R]})$. Similarly, mapping $({{T}}(k),\mu^{(k)},(B^{(k)}_t)_{t\in[0,R]})$ into $l^1$ (from the vertices $(\zeta_i)_{i=1}^{k}$, also defined in the proof of Lemma \ref{subtreelemma}) yields $(\tilde{\mathcal{T}}(k),\tilde{\mu}^{(k)},(\tilde{B}^{(k)}_t)_{t\in[0,1]})$, where $\tilde{B}^{(k)}:=\psi(B^{(k)})$, as in the proof of Proposition \ref{prop1}. Thus the result is a simple consequence of Lemmas \ref{subtreelemma} and \ref{jnkbk}.
\end{proof}

We now consider the convergence of the local times of $\breve{J}^{n,k}$, although before arriving at this result, we must prove a few preparatory lemmas. We will denote the occupation times and local times of $\breve{J}^{n,k}$ by $(\breve{\ell}^{n,k}_m(x))_{m\geq 0,\:x\in V(\breve{\underline{T}}_n(k))}$ and
$(\breve{L}^{n,k}_m(x))_{m\geq 0,\:x\in V(\breve{\underline{T}}_n(k))}$, and define them analogously to (\ref{occtime}) and (\ref{loctime}) respectively. We extend the domains of these processes to the whole of $\mathbb{R}_+\times\breve{\underline{T}}_n(k)$ by linear interpolation, first in space, and then in time. Let us start by proving a simple tail estimate on the occupation times of the jump process.

{\lem \label{occdenstail} Fix $R\in(0,\infty)$ and $k\in\mathbb{N}$. Suppose that Assumption 1 holds, then there exist constants $c_1$, $c_2$ such that, for every $n$,
\[\sup_{x\in V(\breve{{\underline{T}}}_n(k))}\mathbf{P}
 \left(\breve{{\ell}}^{n,k}_{Rn\Lambda_n^{(k)}}(x)\geq n^{1/2}t\right)\leq c_1e^{-c_2t},\hspace{20pt}\forall t\geq 0.\]}
\begin{proof} Using the convergence of the trees proved in Lemma \ref{subtreelemma}, we have, for large $n$, $\min_{i}|e_i^{n,k}|\geq \frac{1}{2}\min_{i}|e_i^{(k)}|$, where the $|e_i^{n,k}|$ are the edge lengths of $\breve{\underline{T}}_n(k)$ and the $|e_i^{(k)}|$ are the edge lengths of ${\underline{T}}_n(k)$. Hence, if $n$ is large enough, for each $x\in V(\breve{\underline{T}}_n(k))$, we can find a line segment of $\breve{\underline{T}}_n(k)$, starting at $x$, which contains no edge endpoints and has length at least $L:=\frac{1}{4}\min_{i}|e_i^{(k)}|$. In particular, it follows from the $n^{-1/2}$ scaling of the trees that this line segment will contain at least $Ln^{1/2}$ vertices in $V(\breve{\underline{T}}_n(k))$. By considering the jump process $\breve{J}^{n,k}$ observed on this line segment, we can use the estimates for the occupation times of a simple random walk on an interval deduced in the appendix (Lemma \ref{Occ}) to obtain an estimate of the appropriate form which holds for large $n$. This is easily extended to $n\in\mathbb{N}$ by suitable choice of $c_1$ and $c_2$, which completes the proof. Note that the two cases considered in Lemma \ref{Occ} cover the possibilities that $x$ is an endpoint of an edge or that it is not.
\end{proof}

We now prove a modulus of continuity result for the local times.

{\lem Fix $R\in(0,\infty)$ and $k\in\mathbb{N}$. Suppose that Assumption 1 holds, then for every $\varepsilon>0$ there exists a constant $c$ such that, for every $n\in\mathbb{N}$, $\delta>0$,
\[\sup_{\buildrel\scriptstyle{x,y\in V(\breve{\underline{{T}}}_n(k)):}\over{d_{\breve{\underline{T}}_n(k)}(x,y)\leq \delta}}\mathbf{P}
 \left(n^{-1/2}\sup_{m\leq Rn\Lambda_n^{(k)}}|\breve{L}^{n,k}_{m}(x)-\breve{L}_m^{n,k}(y)|\geq\varepsilon\right)\leq c\delta^2.\]}
\begin{proof} The argument follows closely the proof of a related estimate in \cite{Borodin}. For brevity we write $R'=R'(n,k)=Rn\Lambda_n^{(k)}$. Fix $x\neq y$ in $V(\breve{\underline{{T}}}_n(k))$ with $d_{\breve{\underline{{T}}}_n(k)}(x,y)\leq \delta$. Conditional on the event where the jump chain $\breve{J}^{n,k}$ hits $x$ before $y$ occurring, we have by a simple calculation
\begin{equation}\label{z}
\sup_{m\leq R'} \left|\breve{L}^{n,k}_{m}(x)-\breve{L}_m^{n,k}(y)-2\sum_{i=1}^{\breve{\ell}_m^{n,k}(x)}\eta_i\right|\leq\sup_{i\leq \breve{\ell}_{R'}^{n,k}(x)+1} 2|\eta_i|,
\end{equation}
where $\eta_i:=N_i\mathrm{deg}_{n,k}(y)^{-1}-\mathrm{deg}_{n,k}(x)^{-1}$. Here, $\mathrm{deg}_{n,k}:=\mathrm{deg}_{\breve{\underline{{T}}}_n(k)}$, and  $N_i$ is the number of visits by $\breve{J}^{n,k}$ to $y$ between the $i$th and $(i+1)$st visits to $x$. Clearly $(\eta_i)_{i\geq 1}$ is an independent identically-distributed family.

Noting that Lemma \ref{occdenstail} allows us to deduce a constant upper bound for the quantity $n^{-1/2}\mathbf{E}(\breve{\ell}_{R'}(x)+1)$ that is uniform in $n$ and $x$, we are able to deduce that
\[{\mathbf{P}\left(\sup_{i\leq \breve{\ell}_{R'}^{n,k}(x)+1} 2|\eta_i|>\varepsilon n^{1/2}\right)}\leq \mathbf{E}\left(\breve{\ell}_{R'}^{n,k}(x)+1\right) \mathbf{P}\left(2|\eta_1|>\varepsilon n^{1/2}\right)\leq  c_1 n^{-3/2}\mathbf{E} |\eta_1|^4,\]
where $c_1$ is a constant that does not depend on $n$, $x$ or $y$. Combining this bound with inequality (\ref{etamoments}) from the appendix implies that
\begin{equation}\label{first}
\mathbf{P}\left(\sup_{i\leq \breve{\ell}_{R'}^{n,k}(x)+1} 2|\eta_i|>\varepsilon n^{1/2}\right)\leq c_2\delta^3,
\end{equation}
where $c_2$ is a constant that does not depend on $n$, $x$ or $y$.

Furthermore, since the sequence $(\sum_{i=0}^m\eta_i)_{m\geq 0}$ is a martingale, we are able to use Doob's martingale norm inequality (see \cite{Kallenberg}, Proposition 7.16, for example) to deduce that
\[\mathbf{P}\left(2\sup_{m\leq R'}\left|\sum_{i=1}^{\breve{\ell}_m^{n,k}(x)}\eta_i\right|>\varepsilon n^{1/2}\right)\leq c_3 n^{-2} \mathbf{E}\left|\sum_{i=1}^{\breve{\ell}_{R'}^{n,k}(x)}\eta_i\right|^4,\]
with $c_3$ not depending on $n$, $x$ or $y$. By replicating the proof of the upper estimate for the corresponding martingale in \cite{Borodin}, using the tail bound of Lemma \ref{occdenstail} and applying the inequality proved in the appendix at (\ref{etamoments}), we are able to bound the right hand side above by $c_4 \delta^2$, uniformly in $n$, $x$ and $y$. Combining this result with (\ref{z}) and the bound at (\ref{first}) yields
\[\mathbf{P}
 \left(\sup_{m\leq R'}|\breve{L}^{n,k}_{m}(x)-\breve{L}_m^{n,k}(y)|\geq\varepsilon n^{1/2}\:\vline\: \min\{m:\breve{J}^{n,k}_m=x\}\leq \min\{m:\breve{J}^{n,k}_m=y\}\right)\leq c_5\delta^2,\]
for every $n$ and $x,y\in V(\breve{\underline{T}}_n(k))$ with $d_{\breve{\underline{T}}_n(k)}(x,y)\leq \delta$. However, if we reverse the role of $x$ and $y$ in the left-hand side, then the same inequality holds, and so we can remove the conditioning to obtain the result.
\end{proof}

We extend this result using a standard maximal inequality.

{\lem \label{smax} Fix $R\in(0,\infty)$ and $k\in\mathbb{N}$. Suppose that Assumption 1 holds, then for every $\varepsilon>0$ there exists a constant $c$ such that, for every $n\in\mathbb{N}$, $\delta>0$,
\begin{equation}\label{fori}
\mathbf{P}
 \left(\sup_{\buildrel\scriptstyle{x,y\in V(\breve{\underline{T}}_n(k)):}\over{d_{\breve{\underline{T}}_n(k)}(x,y)\leq \delta}}n^{-1/2}\sup_{m\leq Rn\Lambda_n^{(k)}}|\breve{L}^{n,k}_{m}(x)-\breve{L}_m^{n,k}(y)|\geq\varepsilon\right)\leq c\delta.
 \end{equation}}
\begin{proof} Let us start by considering a particular edge, $e_i^{k}$ say, of $\underline{T}(k)$. Define $e_i^{n,k}$ to be the corresponding edge in the graph $\breve{\underline{T}}_n(k)$ when the homeomorphism $\Upsilon_{\breve{\underline{T}}_n(k),{\underline{T}}(k)}$ from $\breve{\underline{T}}_n(k)$ to $\underline{T}(k)$ is defined. The set of graph vertices embedded in this edge is $V(e^{n,k}_i):=V(\breve{\underline{T}}_n(k))\cap e_i^{n,k}$. Since an edge of $\breve{\underline{{T}}}_n(k)$ is isomorphic to a Euclidean line-segment, the estimate proved in the previous lemma can be extended by an application of \cite{Bill2}, Theorem 10.3, (or, more precisely, the simple extension of this result that is alluded to in \cite{Bill2}, Problem 10.1), to deduce that
\begin{equation}\label{dok}
\sup_{x \in V(e_i^{n,k})}\mathbf{P} \left(n^{-1/2}\sup_{\buildrel\scriptstyle{y\in V(e_i^{n,k}):}\over{d_{\breve{\underline{T}}_n(k)}(x,y)\leq \delta}}\sup_{m\leq Rn\Lambda_n^{(k)}}|\breve{L}^{n,k}_{m}(x)-\breve{L}_m^{n,k}(y)|\geq\varepsilon\right)\leq c_1\delta^2
\end{equation}
uniformly in $n$ and $\delta$, for some constant $c_1$. Now, since the number of edges of $\breve{\underline{T}}_n(k)$ is bounded uniformly in $n$ for each $k$, there is no problem in replacing the set $V(e_i^{n,k})$ by $V(\breve{\underline{T}}_n(k))$ in the above expression (increasing $c_1$ if necessary).

To complete the proof note that under Assumption 1, for each $n$ and $\delta$, we can choose a $\delta$-net, $A_n^\delta$ say, of $V(\breve{\underline{T}}_n(k))$, such that the quantity $\delta\#A_n^\delta$ is bounded uniformly in $n$ and $\delta$. Applying this fact and the bound at (\ref{dok}) (extended to the whole of $V(\breve{\underline{T}}_n(k))$), it is elementary to check that the left-hand side of (\ref{fori}) is bounded above by
\[\sum_{x\in A_n^\delta}2\mathbf{P} \left(n^{-1/2}\sup_{\buildrel\scriptstyle{y\in V(\breve{\underline{T}}_n(k)):}\over{d_{{\breve{\underline{{T}}}_n(k)}}(x,y)\leq 2\delta }}\sup_{m\leq Rn\Lambda_n^{(k)}}|\breve{L}^{n,k}_{m}(x)-\breve{L}_m^{n,k}(y)|\geq\varepsilon\right)\leq c_2\delta,\]
uniformly in $n$ and $\delta$, which completes the proof.
\end{proof}

We now show that the rescaled local times of the jump-chain $\breve{J}^{n,k}$ are close to those of the Brownian motion $B^{n,k}$ on $(\breve{\underline{T}}_n(k),\lambda_n^{(k)})$. The existence and continuity of the local times of $B^{n,k}$, which we will denote by $\bar{L}^{n,k}$, may be proved by repeating the argument of Lemma \ref{localtimekcont}. The following argument is essentially the same as that used in \cite{Revesz}, Lemma 7, to demonstrate convergence of the local times of the simple random walk on $\mathbb{Z}$.

{\lem \label{iop} Fix $R\in(0,\infty)$ and $k\in\mathbb{N}$. Suppose that Assumption 1 holds and the processes $\breve{J}^{n,k}$ and $B^{n,k}$ are coupled as in the proof of Lemma \ref{jnkbk}, then for every $\varepsilon>0$,
\[\lim_{n\rightarrow\infty}\sup_{x\in V(\breve{\underline{T}}_n(k))}\mathbf{P}\left(\sup_{m\leq Rn\Lambda_n^{(k)}}|n^{-1/2}\breve{L}_m^{n,k}(x)-\bar{L}_{h^{n,k}(m)}^{n,k}(x)|>\varepsilon\right)=0.\]}
\begin{proof} Fix $x\in V(\breve{\underline{T}}_n(k))$. Denote by $(\varsigma_i)_{i\geq 1}$ the hitting times of $x$ by $\breve{J}^{n,k}$, and define $\eta_i:=\bar{L}_{h^{n,k}(\varsigma_i+1)}^{n,k}(x)-\bar{L}_{h^{n,k}(\varsigma_i)}^{n,k}(x)$, where $(h^{n,k}(m))_{m\geq 0}$ are the hitting times defined at (\ref{hitting}). It is straightforward to deduce from the definition of $B^{n,k}$ and the standard scaling properties of one-dimensional Brownian local times that $(\eta_i)_{i\geq 1}$ is an independent, identically-distributed sequence of random variables, each distributed as ${2Z}/{n^{1/2}\mathrm{deg}_{n,k}(x)}$, where $Z$ represents the local time at zero of a standard Brownian motion in $\mathbb{R}^1$, started from zero, evaluated at the hitting time of $\{\pm 1\}$, and $\mathrm{deg}_{n,k}:=\mathrm{deg}_{\breve{\underline{T}}_n(k)}$. The explicit distribution of $Z$ is known as a result of a Ray-Knight theorem (see \cite{Kallenberg}, Theorem 22.17, for example). In particular, $Z$ has finite positive moments of all orders and mean 1. Thus, for $c_1>0$, if we write $R'=R'(n,k)=Rn\Lambda_n^{(k)}$,
\begin{eqnarray}
\label{porb}\lefteqn{\mathbf{P}\left(\sup_{m\leq R'} \left|\eta_1+\dots+\eta_{\breve{\ell}^{n,k}_m(x)}-n^{-1/2}\breve{L}_m^{n,k}(x)\right|>\varepsilon\right)}\\
\nonumber&=& \sum_{m=0}^{\lfloor R'\rfloor} \mathbf{P}\left(\sup_{l\leq m}\left|\eta_1+\dots+\eta_{l}-\frac{2l}{n^{1/2}\mathrm{deg}_{n,k}(x)}\right|>\varepsilon\right)\mathbf{P}\left(\breve{\ell}_{R'}^{n,k}(x)=m\right)\\
\nonumber&\leq & c_2n^{-1}(\ln n)^{2} \sum_{m=0}^{\lfloor c_1n^{1/2} \ln n\rfloor }\mathbf{E}((Z-1)^4) +\mathbf{P}\left(\breve{\ell}_{R'}^{n,k}(x)\geq\lfloor c_1n^{1/2} \ln n\rfloor \right),
\end{eqnarray}
where we have again applied standard martingale inequalities (see \cite{Kallenberg}, Lemma 4.15 and \cite{Bill}, Theorem 6.1, for example) to deduce the inequality. Applying Lemma \ref{occdenstail} and choosing $c_1$ suitably large, we are able to obtain from this an upper bound of the form $c_3 n^{-1/2}(\ln n)^3$ that holds for all $n\geq 2$, uniformly in $x\in V(\breve{\underline{T}}_n(k))$, for the probability at (\ref{porb}).

Observe now that if $\breve{J}_m^{n,k}=x$, then $\eta_1+\dots+\eta_{{\breve{\ell}}^{n,k}_m(x)}=\bar{L}^{n,k}_{h^{n,k}(m+1)}(x)$, otherwise the sum is equal to $\bar{L}^{n,k}_{h^{n,k}(m)}(x)$. Hence
\begin{eqnarray*}
\lefteqn{ \mathbf{P}\left(\sup_{m\leq {R'}}\left|\eta_1+\dots+\eta_{\breve{\ell}^{n,k}_m(x)}-\bar{L}^{n,k}_{h^{n,k}(m)}(x)\right|>\varepsilon\right)}\\
&\leq&{ \mathbf{P}\left(\sup_{m\leq {R'}}\left|\bar{L}^{n,k}_{h^{n,k}(m+1)}(x)-\bar{L}^{n,k}_{h^{n,k}(m)}(x)\right|>\varepsilon\right)}\\
&\leq&\varepsilon^{-4}\sum_{m=0}^{\lfloor R'\rfloor}\mathbf{E}\left(\left|\bar{L}^{n,k}_{h^{n,k}(m+1)}(x)-\bar{L}^{n,k}_{h^{n,k}(m)}(x)\right|^4\right)\\
&=&c_4n^{-1}\mathbf{E}Z^4,
\end{eqnarray*}
uniformly in $x\in V(\breve{\underline{T}}_n(k))$, and since the expectation is finite, this bound converges to zero. The lemma follows.
\end{proof}

We can now combine the estimates of the previous two lemmas to demonstrate that the rescaled local times of $\breve{J}^{n,k}$ converge uniformly to the local times $L^{(k)}$ of $B^{(k)}$. Recall that the domains of the local times $\breve{L}^{n,k}$ are extended to $\mathbb{R}_+\times \breve{\underline{T}}_n(k)$ by linear interpolation.

{\lem \label{localtimeconv} Fix $R\in(0,\infty)$ and $k\in\mathbb{N}$. Suppose that Assumption 1 holds and the processes $\breve{J}^{n,k}$ and $B^{(k)}$ are coupled as in the proof of Lemma \ref{jnkbk}, then for every $\varepsilon>0$,
\[\lim_{n\rightarrow \infty}\mathbf{P}\left(\sup_{x\in{\underline{T}}(k)}\sup_{t\in[0,R]}\left|L^{(k)}_{t}(x)-n^{-1/2}\breve{L}^{n,k}_{tn\Lambda_n^{(k)}}(\Upsilon_{\underline{T}(k),\breve{\underline{T}}_n(k)}(x))\right|\geq\varepsilon\right)=0.\]
}
\begin{proof} In addition to the assumptions of the lemma, suppose also that the processes $B^{n,k}$ and  $B^{(k)}$ are coupled in the way that was used in the proof of Lemma \ref{jnkbk}. As well as the convergence of processes that was described at (\ref{one}), it is possible to show that $\mathbf{P}$-a.s.,
\begin{equation}\label{as}
\lim_{n\rightarrow \infty}\sup_{x\in\underline{T}(k)}\sup_{t\in[0,R]}\left|L^{(k)}_{t}(x)-\bar{L}^{n,k}_t(\Upsilon_{\underline{T}(k),\breve{\underline{T}}_n(k)}(x))\right|=0,
\end{equation}
by first deducing a time-change representation of $\bar{L}^{n,k}\circ\Upsilon_{\underline{T}(k),\breve{\underline{T}}_n(k)}$ in terms of $L^{(k)}$, similar to Lemma \ref{loctimerep}, and then demonstrating that the relevant time-change additive functional converges uniformly in the same way as in Proposition \ref{tildeakt}. This allows the problem to be reduced to showing that
\begin{equation}\label{limit}
\lim_{n\rightarrow\infty}\mathbf{P}\left(\sup_{x\in \breve{\underline{T}}_n(k)}\sup_{t\in[0,R]}\left|n^{-1/2}{\breve{L}}_{tn\Lambda_n^{(k)}}^{n,k}(x)-\bar{L}^{n,k}_{t}(x)\right|\geq\varepsilon\right)=0.
\end{equation}

Now, for each $n$ and $\delta$, we can choose a $\delta$-net, $A_n^\delta$ say, of $\breve{\underline{T}}_n(k)$, consisting of vertices in $V(\breve{\underline{T}}_n(k))$ and such that the quantity $\delta\#A_n^\delta$ is bounded uniformly in $n$ and $\delta$. Using these nets, we can deduce that the probability in the left-hand side of (\ref{limit}) is bounded above by
\begin{eqnarray}
\nonumber\lefteqn{\mathbf{P}\left(\sup_{\buildrel{\scriptstyle{x,y \in \breve{\underline{T}}_n(k):}}\over{d_{\breve{\underline{T}}_n(k)}(x,y)\leq \delta}} \sup_{t\in[0,R]}\left|\bar{{L}}_{t}^{n,k}(x)-\bar{{L}}^{n,k}_{t}(y)\right|\geq\varepsilon/3\right)}\\
\nonumber&+&\sum_{x\in A_n^\delta}\sup_{x\in V(\breve{\underline{T}}_n(k))}\mathbf{P}\left(\sup_{t\in[0,R]}\left|n^{-1/2}\breve{L}_{tn\Lambda_n^{(k)}}^{n,k}(x)-\bar{L}^{n,k}_{t}(x)\right|\geq\varepsilon/3\right)\\
\label{triple}&+&\mathbf{P}\left(\sup_{t\in[0,R]}\sup_{\buildrel{\scriptstyle{x,y \in \breve{\underline{T}}_n(k):}}\over{d_{\breve{\underline{T}}_n(k)}(x,y)\leq \delta}}n^{-1/2}\left|{\breve{L}}^{n,k}_{tn\Lambda_n^{(k)}}(x)-{\breve{L}}^{n,k}_{tn\Lambda_n^{(k)}}(y)\right|\geq\varepsilon/3\right).
\end{eqnarray}
The final term is bounded by $c_1\delta$ uniformly in $n$ by Lemma \ref{smax} (since $\breve{L}^{n,k}$ is extended at each time by linear interpolation over space, there is no problem in extending the result proved there by replacing $V(\breve{\underline{T}}_n(k))$ by $\breve{\underline{T}}_n(k)$). The result at (\ref{as}) implies that the $\limsup$ as $n\rightarrow\infty$ of the first term is bounded above by
\[{\mathbf{P}\left(\sup_{t\in[0,R]}\sup_{\buildrel{\scriptstyle{x,y \in {\underline{T}}(k):}}\over{d_{\underline{T}(k)}(x,y)\leq \delta}} \left|{{L}_t^{(k)}}(x)-{{L}}^{(k)}_{t}(y)\right|\geq\varepsilon/3\right)},\]
and, by choosing $\delta$ appropriately, we can make this probability arbitrarily small since the local times $L^{(k)}$ are jointly continuous in $t$ and $x$. Thus to complete the proof it will suffice to show that the second term of (\ref{triple}) converges to zero for each fixed $\delta$. This is a straightforward consequence of Lemma \ref{iop}, the convergence of local times stated at (\ref{as}), and the strong limit law that was proved for the hitting times $h^{n,k}$ at (\ref{hittingconv}).
\end{proof}

Piecing together the convergence results for trees, measures, jump processes and local times that we have already proved, we obtain the following.

{\propn \label{corr} Fix $R\in(0,\infty)$ and $k\in\mathbb{N}$. Suppose that Assumption 1 holds, then
\[\left(\breve{{T}}_n(k),\breve{\mu}_n^{(k)},\left(\breve{J}^{n,k}_{tn\Lambda_n^{(k)}}\right)_{t\in[0,R]},\left(n^{-1/2}\breve{L}^{n,k}_{tn\Lambda_n^{(k)}}(x)\right)_{t\in[0,R],\:x\in\breve{\underline{{T}}}_n(k)}\right)\]
converges in distribution as $n\rightarrow \infty$ to
\[\left({{T}}(k),{\mu}^{(k)},\left(B^{(k)}_t\right)_{t\in[0,R]},\left({L}^{(k)}_{t}(x)\right)_{t\in[0,R],\:x\in\underline{T}(k)}\right),\]
with respect to the metric $d=(d_1+d_2+d_3+d_4)\wedge 1$.}
\bigskip

From this result, we can deduce the convergence of additive functionals. The definitions of $\hat{A}^{n,k}$ and $\hat{A}^{(k)}$ should be recalled from (\ref{tildeankdef}) and (\ref{tildeakdef}) respectively.

{\cor \label{tildeankak} Fix $R\in(0,\infty)$ and $k\in\mathbb{N}$. Suppose that Assumption 1 holds, then
\[\left(n^{-3/2}\hat{A}^{n,k}_{tn\Lambda_n^{(k)}}\right)_{t\in[0,R]}\Rightarrow \left(\hat{A}^{(k)}_t\right)_{t\in[0,R]}\]
in $C([0,R],\mathbb{R}_+)$.}
\begin{proof} By construction, it will suffice to prove the result when $(\hat{A}^{n,k}_t)_{t\geq 0}$ is replaced by, for $t\geq 0$,
\[\left(n\int_{\breve{\underline{T}}_n(k)}\breve{L}^{n,k}_{(t-1)\vee 0}(x)\breve{\mu}_n^{(k)}(dx)\right)_{t\geq0}.\]
Given Proposition \ref{corr}, this is a straightforward application of the continuous mapping theorem (see \cite{Kallenberg}, Theorem 4.27, for example).
\end{proof}

\section{Tightness for additive functionals}

We now analyse the simple random walks on graph trees in order to obtain a tightness result for the processes $A^{n,k}$ and $\hat{A}^{n,k}$, the definitions of which should be recalled from (\ref{ankdef}) and (\ref{tildeankdef}) respectively. We assume throughout this section that $(\mathcal{T}_n)_{n\geq 1}$ and $(\mathcal{T}_n(k))_{n\geq 1,\: k\geq1}$ are given, and are built from a sequence $\{(w_n,u^n)\}_{n\geq 1}$ that satisfies Assumption 1.

The proof of our main result, Proposition \ref{anktildeank}, is a modification of the argument used by Kesten in \cite{Kestenunpub}, Proposition (4.52), and involves applying some simple random walk estimates for graph trees that are proved in the appendix. In particular, denote the expected holding times of the process $X^{n,k}$ by $\alpha^{n,k}(x):=\mathbf{E}(A_{m+1}^{n,k}-A_m^{n,k}|J^{n,k}_m=x)$, for $x$ a vertex of $\mathcal{T}_n(k)$ and some $m\geq d_{\mathcal{T}_n}(\rho,x)$. Note that the time-homogeneity of the simple random walk means that $\alpha^{n,k}(x)$ is well-defined. By Lemma \ref{alphabeta}, we have the following exact expression for this quantity
\begin{equation}\label{alphaform}
\alpha^{n,k}(x)=\frac{1}{\mathrm{deg}_{n,k}(x)}\left(2n\mu_n^{(k)}(\{x\})-2+\mathrm{deg}_{n,k}(x)\right),
\end{equation}
where we use the notation introduced in Section \ref{discreteproc}, $\mathrm{deg}_{n,k}=\mathrm{deg}_{\mathcal{T}_n(k)}$. We will also consider the expected square value, $\beta^{n,k}(x):=\mathbf{E}((A_{m+1}^{n,k}-A_m^{n,k})^2|J^{n,k}_m=x)$; the bound of Lemma \ref{alphabeta} giving us that
\begin{equation}\label{betaform}
\beta^{n,k}(x)\leq 36n^2(\deg_{n,k}(x)+\Delta_n^{(k)})\frac{\mu_n^{(k)}(\{x\})^2}{\mathrm{deg}_{n,k}(x)}.
\end{equation}

Before continuing, for want of a suitable reference we state a simple lemma, which may be proved using elementary probability theory. See \cite{Kallenberg}, Exercise 6.11, for a closely related result.

{\lem \label{condconvlem} Let $(Z^{n,k})_{n,k\geq 1}$ be a collection of random variables and $(\mathcal{F}^{n,k})_{n,k\geq 1}$ a collection of $\sigma$-algebras on the probability space with probability measure $\mathbf{P}$. If, for every $\varepsilon>0$,
\[\lim_{k\rightarrow\infty}\limsup_{n\rightarrow\infty}\mathbf{P}\left(\mathbf{E}(Z^{n,k}|\mathcal{F}^{n,k})>\varepsilon\right)=0,\]
then, for every $\varepsilon>0$,
\[\lim_{k\rightarrow\infty}\limsup_{n\rightarrow\infty}\mathbf{P}\left(Z^{n,k}>\varepsilon\right)=0.\]}

In the time-scaling of the following result, it will be useful to include the quantity $\Lambda_n^{(k)}$, which was introduced at (\ref{length}). Note that, under Assumption 1, the limit as $n\rightarrow\infty$ of $\Lambda_n^{(k)}$ exists and is finite for each fixed $k\in \mathbb{N}$.

{\propn \label{anktildeank} Fix $R\in(0,\infty)$ and $\varepsilon>0$. If Assumption 1 holds, then
\[\lim_{k\rightarrow\infty}\limsup_{n\rightarrow\infty}\mathbf{P}\left(n^{-3/2}\sup_{m\leq Rn\Lambda_n^{(k)}}|{A}^{n,k}_{m}-\hat{A}^{n,k}_{m}|>\varepsilon\right)=0.\]}
\begin{proof} Let $m\in\mathbb{N}$. By definition, we have that
\begin{eqnarray}
|{A}^{n,k}_{m}-\hat{A}^{n,k}_{m}|&=&\left|\sum_{l=0}^{m-1}\left(A_{l+1}^{n,k}-A_l^{n,k}-\frac{2n\mu_n^{(k)}(\{J_l^{n,k}\})}{\mathrm{deg}_{n,k}(J_l^{n,k})}\right)\right|\nonumber\\
&\leq& \left|\sum_{l=0}^{m-1}\left(A_{l+1}^{n,k}-A_l^{n,k}-\alpha^{n,k}(J_l^{n,k})\right)\right|\nonumber\\
&&\hspace{40pt}+\sum_{l=0}^{m-1}\left|\alpha^{n,k}(J^{n,k}_l)-\frac{2n\mu_n^{(k)}(\{J_l^{n,k}\})}{\mathrm{deg}_{n,k}(J_l^{n,k})}\right|.\label{sumyo}
\end{eqnarray}
We shall consider these two terms separately, starting with the second summand. First, we use the formula at (\ref{alphaform}) in place of $\alpha^{n,k}$ to deduce that
\[\sup_{m\leq Rn\Lambda_n^{(k)}}\sum_{l=0}^{m-1}\left|\alpha^{n,k}(J^{n,k}_l)-\frac{2n\mu_n^{(k)}(\{J_l^{n,k}\})}{\mathrm{deg}_{n,k}(J_l^{n,k})}\right|= \sum_{l=0}^{\lfloor Rn\Lambda_n^{(k)}\rfloor-1} \frac{|2-\mathrm{deg}_{n,k}(J_l^{n,k})|}{\mathrm{deg}_{n,k}(J_l^{n,k})}\leq  Rn\Lambda_n^{(k)}.\]
Hence, when multiplied by $n^{-3/2}$, as $n\rightarrow\infty$, the second term of (\ref{sumyo}) converges to zero uniformly in $m\leq Rn\Lambda_n^{(k)}$, $\mathbf{P}$-a.s.

We now deal with the first summand of (\ref{sumyo}). Since, conditional on knowing $J^{n,k}$, the expected value of $A_{l+1}^{n,k}-A_l^{n,k}$ is precisely $\alpha^{n,k}(J_l^{n,k})$, we can use Kolmogorov's maximum inequality (see \cite{Kallenberg}, Lemma 4.15) to deduce that, for $\varepsilon>0$,
\begin{eqnarray}
\nonumber\lefteqn{\mathbf{P}\left(n^{-3/2}\sup_{m\leq Rn\Lambda_n^{(k)}}\left|\sum_{l=0}^{m-1}\left(A_{l+1}^{n,k}-A_l^{n,k}-\alpha^{n,k}(J_l^{n,k})\right)\right|>\varepsilon\:\vline\: J^{n,k}\right)}\\
\nonumber &\leq &\frac{1}{n^{3}\varepsilon^2}\sum_{l=0}^{\lfloor Rn\Lambda_n^{(k)}\rfloor-1}\beta^{n,k}(J_l^{n,k})\\
\nonumber &\leq&\frac{36}{n\varepsilon^2}\sum_{l=0}^{\lfloor Rn\Lambda_n^{(k)}\rfloor-1}(\mathrm{deg}_{n,k}(J_l^{n,k})+\Delta_n^{(k)})\frac{\mu_n^{(k)}(\{J_l^{n,k}\})}{\mathrm{deg}_{n,k}(J_l^{n,k})}\hspace{140pt}\\
\label{ub} &\leq&\frac{18}{n^2\varepsilon^2}\hat{A}^{n,k}_{\lfloor Rn\Lambda_n^{(k)}\rfloor}\left(\max_{x\in\mathcal{T}_n(k)}\mathrm{deg}_{n,k}(x)+\Delta_n^{(k)}\right),
\end{eqnarray}
where we have used the bound at (\ref{betaform}) for the second inequality, and we have also dropped a power of $\mu_n^{(k)}(\{x\})$, which is allowed because $\mu_n^{(k)}(\{x\})\leq 1$. The final inequality follows simply from the definition of $\hat{A}^{n,k}$. For $\delta>0$, we have
\begin{eqnarray}
\nonumber \lefteqn{\lim_{k\rightarrow\infty}\limsup_{n\rightarrow\infty}\mathbf{P}\left(n^{-2}\hat{A}^{n,k}_{\lfloor Rn\Lambda_n^{(k)}\rfloor}\left(\max_{x\in\mathcal{T}_n(k)}\mathrm{deg}_{n,k}(x)+\Delta_n^{(k)}\right)>\delta\right)}\\
\label{twot} &\leq&\lim_{k\rightarrow\infty}\limsup_{n\rightarrow\infty}\left[\frac{R+1}{\delta n^{1/2}}\left(\max_{x\in\mathcal{T}_n(k)}\mathrm{deg}_{n,k}(x)+\Delta_n^{(k)}\right)+\mathbf{P}\left(n^{-3/2}\hat{A}^{n,k}_{\lfloor Rn\Lambda_n^{(k)}\rfloor}>R+1\right)\right]\hspace{15pt}
\end{eqnarray}
Now it is a simple consequence of Lemma \ref{subtreelemma} that $\max_{x\in\mathcal{T}_n(k)}\mathrm{deg}_{n,k}(x)$ is bounded uniformly in $n$. Combined with Lemma \ref{equiv}, this implies that
\[\lim_{k\rightarrow\infty}\limsup_{n\rightarrow\infty}n^{-1/2}\left(\max_{x\in\mathcal{T}_n(k)}\mathrm{deg}_{n,k}(x)+\Delta_n^{(k)}\right)=0,\]
which deals with the first of the terms of (\ref{twot}). To show the second term is also zero, we apply the distributional convergence results of Proposition \ref{tildeakt} and Corollary \ref{tildeankak}. Hence Lemma \ref{condconvlem} allows us to deduce from the upper bound at (\ref{ub}) that
\[\lim_{k\rightarrow\infty}\limsup_{n\rightarrow\infty}\mathbf{P}\left(n^{-3/2}\sup_{m\leq Rn\Lambda_n^{(k)}}\left|\sum_{l=0}^{m-1}\left(A_{l+1}^{n,k}-A_l^{n,k}-\alpha^{n,k}(J_l^{n,k})\right)\right|>\varepsilon\right)=0.\]
By recalling the bound for $|A_m^{n,k}-\hat{A}_m^{n,k}|$ from (\ref{sumyo}), and applying the limit results that we have proved for each of the summands, it is straightforward to deduce the desired result.
\end{proof}

In conjunction with the convergence results we have already proved for $\hat{A}^{(k)}$ and $\hat{A}^{n,k}$ in Proposition \ref{tildeakt} and Corollary \ref{tildeankak} respectively, from the above proposition we are able to deduce a concrete description of the growth of $A^{n,k}$ as $n$ and then $k$ get large. We assume that $A^{n,k}$ is extended to a continuous time process by linear interpolation.

{\cor\label{agrowth} Fix $R\in(0,\infty)$ and $\varepsilon>0$. If Assumption 1 holds, then
\[\lim_{k\rightarrow\infty}\limsup_{n\rightarrow\infty}\mathbf{P}\left(\sup_{t\leq R}|n^{-3/2}{A}^{n,k}_{tn\Lambda_n^{(k)}}-t|>\varepsilon\right)=0.\]}

\section{Tightness of discrete processes}

As in the previous section, we assume that $(\mathcal{T}_n)_{n\geq 1}$ and $(\mathcal{T}_n(k))_{n\geq 1,\: k\geq1}$ are given, and are constructed from a sequence $\{(w_n,u^n)\}_{n\geq 1}$ that satisfies Assumption 1. Consequently we can define the isometric embedding $\psi_n:\mathcal{T}_n\rightarrow l^1$ as at the end of Section \ref{cbptree}. We shall denote the $l^1$-embedded versions of $X^n$, $X^{n,k}$ and $J^{n,k}$ by $\tilde{X}^{n}$, $\tilde{X}^{n,k}$ and $\tilde{J}^{n,k}$ respectively, and extend the definitions of these discrete time processes to continuous time by linear interpolation. The main result of this section is obtained in Corollary \ref{tightproc}, which demonstrates a tightness result for $\tilde{X}^n$ and $\tilde{J}^{n,k}$ when these processes are rescaled appropriately. We start by proving a lemma which provides a modulus of continuity result for the jump processes.

{\lem \label{jmodcont} Fix $R\in(0,\infty)$ and $\varepsilon>0$. If Assumption 1 holds, then
\[\lim_{\delta\rightarrow 0}\limsup_{k\rightarrow\infty}\limsup_{n\rightarrow\infty}\mathbf{P}\left(n^{-1/2} \sup_{s,t\leq R: \: |s-t|\leq\delta} \|\tilde{J}_{sn\Lambda_n^{(k)}}^{n,k}-\tilde{J}_{tn\Lambda_n^{(k)}}^{n,k}\|\geq \varepsilon\right)=0.\]}
\begin{proof} By the convergence results of Propositions \ref{prop1} and \ref{jnkbkl1}, it is sufficient to show that
\[\lim_{\delta\rightarrow 0} \mathbf{P}\left( \sup_{s,t\leq R: \: |s-t|\leq\delta}\|\tilde{X}_s-\tilde{X}_t\|\geq\varepsilon\right)=0,\]
where $\tilde{X}$ is the $l^1$-embedded version of $X$ defined in the proof of Proposition \ref{prop1}. This is a simple consequence of the fact that $\tilde{X}$ is continuous, $\mathbf{P}$-a.s.
\end{proof}

We now have enough information to demonstrate a tightness result for $\tilde{X}^{n,k}$ and $\tilde{J}^{n,k}$.

{\propn  Fix $\varepsilon>0$. If Assumption 1 holds, then
\[\lim_{k\rightarrow\infty}\limsup_{n\rightarrow \infty}\mathbf{P}\left(n^{-1/2}\sup_{t\in[0,1]}\|\tilde{X}^{n,k}_{tn^{3/2}}-\tilde{J}^{n,k}_{tn\Lambda_n^{(k)}}\|>\varepsilon\right)=0.\]}
\begin{proof} Fix $\varepsilon,\eta>0$. By the modulus of continuity result of Lemma \ref{jmodcont}, we can choose $\delta\in(0,1)$ such that
\begin{equation}\label{e2}
\limsup_{k\rightarrow\infty}\limsup_{n\rightarrow\infty}\mathbf{P}\left(n^{-1/2} \sup_{s,t\leq 3: \: |s-t|\leq\delta} \|\tilde{J}_{sn\Lambda_n^{(k)}}^{n,k}-\tilde{J}_{tn\Lambda_n^{(k)}}^{n,k}\|\geq \varepsilon\right)\leq \eta.
\end{equation}
Set $E_1(n,k):=\{n^{-1/2} \sup_{s,t\leq 3: \: |s-t|\leq\delta} \|\tilde{J}_{sn\Lambda_n^{(k)}}^{n,k}-\tilde{J}_{tn\Lambda_n^{(k)}}^{n,k}\|< \varepsilon\}$. Also define
\[E_2(n,k):=\left\{\sup_{t\leq 2}|n^{-3/2}{A}^{n,k}_{tn\Lambda_n^{(k)}}-t|<\delta\right\},\]
Note that by (\ref{e2}) and Corollary \ref{agrowth},
\begin{equation}\label{e12345}
\limsup_{k\rightarrow\infty}\limsup_{n\rightarrow\infty}\mathbf{P}\left(\bigcup_{i=1}^2 E_i(n,k)^c\right)\leq \eta.
\end{equation}

Assume for the next part of the argument that $\cap_{i=1}^2 E_i(n,k)$ holds, and note that on $E_2(n,k)$ we have, for $t\in[\delta,1]$,
\[A_{(t-\delta)n\Lambda_n^{(k)}}^{n,k}<tn^{3/2}<A_{(t+\delta)n\Lambda_n^{(k)}}^{n,k}.\]
Now recall the definition of $\tau^{n,k}$ from (\ref{taudef}), and note that, because $A_t^{n,k}$ is strictly increasing and linear between integer times, then if we extend the definition of $\tau^{n,k}$ to continuous time by linear interpolation, then $\tau^{n,k}$ satisfies $\tau^{n,k}(t):=\max\{s:\:n^{-3/2}A_s^{n,k}\leq t\}$ for $t\geq 0$.  As a simple consequence of this and the above pair of inequalities, it must be the case that $|{\tau}^{n,k}(tn^{3/2})- tn\Lambda_n^{(k)}|\leq \delta n \Lambda_n^{(k)}$, for $t\in[0,1]$. On $E_1(n,k)$, we have a bound for the modulus of continuity of the jump process $\tilde{J}^{n,k}$, and using the previous inequality, it is possible to deduce from this that
\[n^{-1/2}\sup_{t\in[0,1]}\|\tilde{J}^{n,k}_{{\tau}^{n,k}(tn^{3/2})}-\tilde{J}^{n,k}_{tn\Lambda_n^{(k)}}\|< \varepsilon.\]
However, after relabeling using (\ref{recover}), we are able to obtain from this that
\[n^{-1/2}\sup_{t\in[0,1]}\|\tilde{X}^{n,k}_{tn^{3/2}}-\tilde{J}^{n,k}_{tn\Lambda_n^{(k)}}\|\leq \varepsilon+n^{-1/2},\]
where the extra $n^{-1/2}$ term arises due to the difference in the linear interpolation procedures used when defining the processes $\tilde{J}^{n,k}_{{\tau}^{n,k}(tn^{3/2})}$ and $\tilde{X}^{n,k}_{tn^{3/2}}$. Thus we reach the conclusion that
\[{\limsup_{k\rightarrow\infty}\limsup_{n\rightarrow \infty}\mathbf{P}\left(\sup_{t\in[0,1]}\|\tilde{X}^{n,k}_{tn^{3/2}}-\tilde{J}^{n,k}_{tn\Lambda_n^{(k)}}\|>2\varepsilon n^{1/2}\right)}
\leq \limsup_{k\rightarrow\infty}\limsup_{n\rightarrow\infty}\mathbf{P}\left(\bigcup_{i=1}^2 E_i(n,k)^c\right)\]
which, by (\ref{e12345}), is bounded above by $\eta$. Since $\eta$ was arbitrary, the proof is complete.
\end{proof}

{\lem \label{lemyo} If Assumption 1 holds, then $\mathbf{P}$-a.s.
\[\lim_{k\rightarrow\infty}\limsup_{n\rightarrow \infty}n^{-1/2}\sup_{t\in[0,1]}\|\tilde{X}^{n}_{tn^{3/2}}-{\tilde{X}}^{n,k}_{tn^{3/2}}\|=0.\]}
\begin{proof} From the definition of the process $X^{n,k}$ as the projection of $X^n$ onto $\mathcal{T}_n(k)$, it is clear that the supremum in the expression is bounded by $\Delta_n^{(k)}$, as defined at (\ref{deltankdef}). Hence the result follows from Lemma \ref{equiv}.
\end{proof}

The two previous results immediately imply the following.

{\cor \label{tightproc} If Assumption 1 holds, then
\[\lim_{k\rightarrow\infty}\limsup_{n\rightarrow \infty}\mathbf{P}\left(n^{-1/2}\sup_{t\in[0,1]}\|\tilde{X}^{n}_{tn^{3/2}}-\tilde{J}^{n,k}_{tn\Lambda_n^{(k)}}\|>\varepsilon\right)=0.\]}

\section{Convergence of quenched law}

All the hard analysis of the proof of Theorem \ref{mainresult} is now complete. However, before proving it, we summarise the tightness result for the laws of the rescaled processes that we will apply. As in the previous section, we use the notation $\tilde{J}^{n,k}=\psi_n(J^{n,k})$ and $\tilde{X}^{n}=\psi_n(X^n)$, where $\psi_n$ is the distance-preserving embedding of vertices of $\mathcal{T}_n$ into $l^1$ described at the end of Section \ref{cbptree}, and these discrete time processes are extended to continuous time by linear interpolation. We also include the corresponding tightness results for sets and measures.

{\propn \label{tite} Suppose Assumption 1 holds. If we denote by $\tilde{\mathbf{Q}}^{\mathcal{T}_n(k)}_\rho$ the law of the process $(n^{-1/2}\tilde{J}^{n,k}_{tn\Lambda_n^{(k)}})_{t\in [0,1]}$, and by $\tilde{\mathbf{Q}}^{\mathcal{T}_n}_\rho$ the law of $(n^{-1/2}\tilde{X}^{n}_{tn^{3/2}})_{t\in [0,1]}$, then
\[\lim_{k\rightarrow\infty}\limsup_{n\rightarrow\infty}d_H^{l^1}(n^{-1/2}\tilde{\mathcal{T}}_n,n^{-1/2}\tilde{\mathcal{T}}_n(k))=0,\]
\[\lim_{k\rightarrow\infty}\limsup_{n\rightarrow\infty}d_P^{l^1}(\tilde{\mu}_n(n^{1/2}\cdot),\tilde{\mu}^{(k)}_n(n^{1/2}\cdot))=0,\]
\[\lim_{k\rightarrow\infty}\limsup_{n\rightarrow\infty}d_P^{C([0,1],l^1)}(\tilde{\mathbf{Q}}^{\mathcal{T}_n}_\rho,\tilde{\mathbf{Q}}^{\mathcal{T}_n(k)}_\rho)=0,\]
where $d_H^{l^1}$ is the Hausdorff metric on $\mathcal{K}(l^1)$, and $d_P^{\cdot}$ is the Prohorov metric on $\mathcal{M}_1(\cdot)$.}
\begin{proof} The first two limits are consequences of Lemma \ref{equiv}, and the definitions of $\mathcal{T}_n(k)$ and $\mu_n^{(k)}$ using the projection operator. The third limit can be deduced from Corollary \ref{tightproc}.
\end{proof}

{\propn \label{bigresult} If Assumption 1 holds, then
\[\Theta_n\left(\tilde{\mathcal{T}}_n,\tilde{\mu}_n,\tilde{\mathbf{P}}^{\mathcal{T}_n}_\rho\right)\rightarrow(\tilde{\mathcal{T}},\tilde{\mu},\tilde{\mathbf{P}}^\mathcal{T}_\rho)\]
in the space $\mathcal{K}(l^1)\times\mathcal{M}_1(l^1)\times \mathcal{M}_1(C([0,1],l^1))$.}
\begin{proof} Elementary analysis can be used to obtain the result from Propositions \ref{prop1}, \ref{jnkbkl1} and \ref{tite}.
\end{proof}

\hspace{10pt}
\begin{proof11} Suppose $n^{-1/2}w_n\rightarrow w\in \mathcal{W}^*$. By definition of $\mathcal{W}^*$, see (\ref{wstar}), there exists a $u\in [0,1]^{\mathbb{N}}$, such that $(w,u)\in\Gamma$. If we now take $u^n=u$ for each $n$, then the sequence $(w_n,u^n)$ satisfies Assumption 1, with the relevant limit being given by $(w,u)$. Thus Theorem \ref{mainresult} follows from Proposition \ref{bigresult}.
\end{proof11}

\section{Measurability and convergence of annealed law}\label{meassec}

Given the quenched limit result of Theorem \ref{mainresult}, there is little to do to establish the annealed limit of Theorem \ref{ann} apart from check the measurability of various objects, and that is the primary aim of this section. Note that in all the discussions of measurability that follow, we assume that the $\sigma$-algebra of the underlying probability space is $\mathbf{P}$-complete (which is no real restriction, as we can easily complete it if it is not already). Furthermore, to avoid confusion we will apply subscripts to objects built from deterministic pairs $(w,u)\in C([0,1],\mathbb{R}_+)\times [0,1]^\mathbb{N}$, as in Section \ref{abstracttree}, in the following way: $\tilde{\mathcal{T}}_{w,u},\tilde{\mu}_{w,u},\tilde{\mathbf{P}}_\rho^{\mathcal{T}_{w,u}},\dots$.

We start by showing that the $l^1$-embedded triple $(\tilde{\mathcal{T}},\tilde{\mu},\tilde{\mathbf{P}}_\rho^{\mathcal{T}}):=(\tilde{\mathcal{T}}_{W,U},\tilde{\mu}_{W,U},\tilde{\mathbf{P}}_\rho^{\mathcal{T}_{W,U}})$ is $(W,U)$-measurable, where $(W,U)$ are the random variables defined at the start of Section \ref{crtsec}. Since we have only defined $(\tilde{\mathcal{T}}_{w,u},\tilde{\mu}_{w,u},\tilde{\mathbf{P}}_\rho^{\mathcal{T}_{w,u}})$ so far for $(w,u)\in \Gamma$, we extend the definition to the entire of $C([0,1],\mathbb{R}_+)\times [0,1]^{\mathbb{N}}$ by setting it to be an arbitrary constant triple on the set $\Gamma^c$. The notation $\tilde{\mathbf{P}}^{\mathcal{T}(k)}_\rho$ should be recalled from (\ref{tpk}).

{\lem \label{measlem1} (a) For each $k\in \mathbb{N}$, the map from $\Gamma\subseteq C([0,1],\mathbb{R}_+)\times [0,1]^\mathbb{N}$ (equipped with the usual subspace topology) to $\mathcal{K}(l^1)\times\mathcal{M}_1(l^1)\times\mathcal{M}_1(C([0,1],l^1))$ that takes the pair $(w,u)$ to $(\tilde{\mathcal{T}}_{w,u}(k),\tilde{\mu}_{w,u}^{(k)},\tilde{\mathbf{P}}_\rho^{\mathcal{T}_{w,u}(k)})$ is continuous.\\
(b) The map $(w,u)\mapsto(\tilde{\mathcal{T}}_{w,u},\tilde{\mu}_{w,u},\tilde{\mathbf{P}}_\rho^{\mathcal{T}_{w,u}})$ defines a measurable function from $\Gamma$ (equipped with the subspace $\sigma$-algebra) to $\mathcal{K}(l^1)\times\mathcal{M}_1(l^1)\times\mathcal{M}_1(C([0,1],l^1))$.\\
(c) The triple $(\tilde{\mathcal{T}},\tilde{\mu},\tilde{\mathbf{P}}_\rho^{\mathcal{T}})$ is $(W,U)$-measurable.}
\begin{proof} Let $(w^n,u^n)\in\Gamma$ be such that $(w^n,u^n)\rightarrow (w,u)\in\Gamma$. By repeating an almost identical argument to Lemma \ref{subtreelemma} (and mapping this result into $l^1$ using the sequential construction), we are able to show that $(\tilde{\mathcal{T}}_{w^n,u^n}(k),\tilde{\mu}_{w^n,u^n}^{(k)})\rightarrow(\tilde{\mathcal{T}}_{w,u}(k),\tilde{\mu}_{w,u}^{(k)})$, which deals with the first two coordinates. The simultaneous convergence of laws in $\mathcal{M}_1(C([0,1],l^1))$ can be proved by following the steps that lead to (\ref{one}), and then mapping into $l^1$. This completes the proof of part (a), which has as a consequence that $(w,u)\mapsto(\tilde{\mathcal{T}}_{w,u}(k),\tilde{\mu}_{w,u}^{(k)},\tilde{\mathbf{P}}_\rho^{\mathcal{T}_{w,u}(k)})$ is measurable on $\Gamma$. Recall from Proposition \ref{prop1} that on $\Gamma$ we have
\[\left(\tilde{\mathcal{T}}_{w,u}(k),\tilde{\mu}_{w,u}^{(k)},\tilde{\mathbf{P}}_\rho^{\mathcal{T}_{w,u}(k)}\right) \rightarrow\left(\tilde{\mathcal{T}}_{w,u},\tilde{\mu}_{w,u},\tilde{\mathbf{P}}_\rho^{\mathcal{T}_{w,u}}\right).\] Since a limit of measurable functions is again measurable, this implies part (b). Finally, applying the fact that $\Gamma$ is a measurable subset of $C([0,1],\mathbb{R}_+)\times [0,1]^\mathbb{N}$ chosen (in Lemma \ref{crtprops}) to satisfy $\mathbf{P}((W,U)\in\Gamma)=1$, part (c) follows easily.
\end{proof}

This result allows us to deduce the existence of a probability measure satisfying (\ref{pdef}). First, denote by $\Omega$ our underlying probability space, so that $(W,U)=(W(\omega),U(\omega))$ and $\tilde{\mathbf{P}}_\rho^\mathcal{T}=\tilde{\mathbf{P}}_\rho^{\mathcal{T}_{W(\omega),U(\omega)}}$, where $\omega\in\Omega$. By part (c) of the above lemma, the collection of laws $(\tilde{\mathbf{P}}^{\mathcal{T}}_\rho)_{\omega\in\Omega}$ can be viewed as a probability kernel from $\Omega$ to $C([0,1],l^1)$, (see \cite{Kallenberg}, Lemma 1.40). Thus we can extend the probability measure $\mathbf{P}$ on $\Omega$ to a probability measure $\hat{\mathbf{P}}$ on $\Omega\times C([0,1],l^1)$ by setting
\begin{equation}\label{extend}
\hat{\mathbf{P}}(d\omega d\tilde{X}):=\mathbf{P}(d\omega)\tilde{\mathbf{P}}^{\mathcal{T}}_\rho(d\tilde{X}),
\end{equation}
for $\omega\in\Omega$, $\tilde{X}\in C([0,1],l^1)$. The above lemma also allows us to deduce that $(\omega,\tilde{X})\mapsto(W(\omega),U(\omega),\tilde{\mathcal{T}}_{W(\omega),U(\omega)},\tilde{\mu}_{W(\omega),U(\omega)},\tilde{X})$ is a measurable function on $\Omega\times C([0,1],l^1)$, and moreover
\begin{eqnarray}
\nonumber\lefteqn{\hat{\mathbf{P}}\left(\tilde{\mathcal{T}}\in A,\:\tilde{\mu}\in B,\:\tilde{X}\in C\right)}\\
\nonumber&=&\int_{\Omega\times C([0,1],l^1)}\mathbf{P}(d\omega)\tilde{\mathbf{P}}^{\mathcal{T}}_\rho(d\tilde{X})\mathbf{1}_{\{\tilde{\mathcal{T}}\in A,\:\tilde{\mu}\in B,\:\tilde{X}\in C\}}\\
\nonumber&=&\int_{\Omega}\mathbf{P}(d\omega)\mathbf{1}_{\{\tilde{\mathcal{T}}\in A,\:\tilde{\mu}\in B\}}\tilde{\mathbf{P}}^{\mathcal{T}}_\rho(C)\\
\label{forty}&=&\int_{C([0,1],\mathbb{R}_+)\times [0,1]^\mathbb{N}} \mathbf{P}((W,U)\in (dw,du))\: \mathbf{1}_{\{\tilde{\mathcal{T}}\in A,\:\tilde{\mu}\in B\}}\tilde{\mathbf{P}}^\mathcal{T}_\rho(C),
\end{eqnarray}
for every measurable $A\subseteq\mathcal{K}(l^1)$, $B\subseteq\mathcal{M}_1(l^1)$, and $C\subseteq C([0,1],l^1)$, where the final equality is obtained by a simple change of variables in the integral. Hence if we define $\mathbb{P}$ to be the law of $(\tilde{\mathcal{T}},\tilde{\mu},\tilde{X})$ under the measure $\hat{\mathbf{P}}$, then $\mathbb{P}$ satisfies (\ref{pdef}). That it is the unique measure to do so is standard (see \cite{Kallenberg}, Lemma 1.17, for example). Finally, that the law of $\tilde{X}$ under the conditional measure $\hat{\mathbf{P}}(\cdot|(W,U))$ is given by $\tilde{\mathbf{P}}_\rho^\mathcal{T}$ is readily deduced from (\ref{forty}).

To prove the corresponding discrete results we can follow similar arguments, and so we will omit the proofs. Henceforth, we suppose that $(\mathcal{T}_n)_{n\geq 1}$ is a sequence of random ordered graph trees whose search-depth functions $(W_n)_{n\geq 1}$ are independent of $U$, and also satisfy the convergence result at (\ref{sdc}). The triple $(\tilde{\mathcal{T}}_n,\tilde{\mu}_n,\tilde{\mathbf{P}}_\rho^{\mathcal{T}_n})\in\mathcal{K}(l^1)\times\mathcal{M}_1(l^1)\times \mathcal{M}_1(C(\mathbb{R}_+,l^1))$ is constructed from the random pair $(W_n,U)$ by following the procedure presented in Sections {\ref{cbptree} and \ref{discreteproc} for deterministic pairs $(w_n,u)$, and analogously to Lemma \ref{measlem1} we have that this construction is measurable. By extending the underlying probability space in a similar fashion to (\ref{extend}), we can also deduce the existence of a probability measure $\mathbb{P}_n$ satisfying (\ref{pndef}).

To complete this section, we prove the annealed limit result of Theorem \ref{ann}. The two versions of the definition (one involving laws, and one involving processes) of the rescaling operator $\Theta_n$ should be recalled from the introduction.

\hspace{20pt}
\begin{proof12} By \cite{Kallenberg}, Theorem 4.29, it is sufficient to demonstrate that $\mathbb{P}_n\circ \Theta_n^{-1}(F)\rightarrow\mathbb{P}(F)$ for any function $F$ of the form $F(K,\nu,f)=F_1(K)F_2(\nu)F_3(f)$, where $F_1\in C_b(\mathcal{K}(l^1))$, $F_2\in C_b(\mathcal{M}_1(l^1))$ and $F_3\in C_b(C([0,1],l^1))$. Thus we fix $F$ to be such a function.

Now by assumption we have that $n^{-1/2}W_n\Rightarrow W$, and so $(n^{-1/2}W_n,U)\Rightarrow (W,U)$. As a result of the separability of $C([0,1],\mathbb{R}_+)\times [0,1]^{\mathbb{N}}$, it is therefore possible to construct $(n^{-1/2}W_n^*,U_n^*)$, which is a version of $(n^{-1/2}W_n,U)$, and $(W^*,U^*)$, which is a version of $(W,U)$, in such a way that $(n^{-1/2}W_n^*,U_n^*)\rightarrow(W^*,U^*)$ almost-surely on some probability space, $\Omega^*$ say, with probability measure $\mathbf{P}^*$. We now suppose that the random triple $(\tilde{\mathcal{T}},\tilde{\mu},\tilde{\mathbf{P}}_\rho^{\mathcal{T}})$ is built from $(W^*,U^*)$ and that the random triple $(\tilde{\mathcal{T}}_n,\tilde{\mu}_n,\tilde{\mathbf{P}}_\rho^{\mathcal{T}_n})$ is built from $(W_n^*,U_n^*)$ for each $n$.

It is easy to check that the random variables $(n^{-1/2}W_n^*,U_n^*)_{n\geq 1}$ satisfy Assumption 1 with the relevant limit being given by $(W^*,U^*)\in \Gamma$, $\mathbf{P}^*$-a.s. As a consequence of this, Proposition \ref{bigresult} implies that $\mathbf{P}^*$-a.s., $\Theta_n(\tilde{\mathcal{T}}_n,\tilde{\mu}_n,\tilde{\mathbf{P}}_\rho^{\mathcal{T}_n})\rightarrow(\tilde{\mathcal{T}},\tilde{\mu},\tilde{\mathbf{P}}_\rho^{\mathcal{T}})$.
Thus we have that, $\mathbf{P}^*$-a.s.,
$F_1(n^{1/2}\tilde{\mathcal{T}_n})\rightarrow F_1(\tilde{\mathcal{T}})$, $F_2(\tilde{\mu}_n(n^{1/2}\cdot))\rightarrow F_2(\tilde{\mu})$ and
\[\int_{\tilde{X}\in C(\mathbb{R}_+,l^1)} \tilde{\mathbf{P}}_\rho^{\mathcal{T}_n}(df)F_3\left((n^{-1/2}\tilde{X}(tn^{3/2}))_{t\in[0,1]}\right)\rightarrow \tilde{\mathbf{P}}_\rho^{\mathcal{T}}(F_3).\]
By applying the dominated convergence theorem (twice), it follows that $\mathbb{P}_n\circ\Theta_n^{-1}(F)\rightarrow \mathbb{P}(F)$, as desired.
\end{proof12}

\appendix
\section*{Appendix}

\section{Index of processes}\label{notation}

The list below provides a summary (in order of first appearance) of the more important random processes that appear in the article.
\bigskip

\begin{center}
\begin{tabular}{clc}
  $X$ & Brownian motion on $(\mathcal{T},\mu)$. & Section \ref{crtproc}\\
  $L$ & Local times of $X$. & Lemma \ref{localtimecont}\\
  $A^{(k)}$ & Time-change additive functional from  $X$ to $B^{(k)}$. & (\ref{akdef})\\
  $\tau^{(k)}$ & Inverse of $A^{(k)}$. &(\ref{taukdef})\\
  $B^{(k)}$ & Brownian motion on $(\mathcal{T}(k),\lambda^{(k)})$. & (\ref{bkdef})\\
  $X^n$ & Simple random walk on $\mathcal{T}_n$. & Section \ref{discreteproc}\\
  $X^{n,k}$ & Projection of $X^n$ onto $\mathcal{T}_n(k)$. & (\ref{xnkdef})\\
  $J^{n,k}$ & Jump process associated with $X^{n,k}$. & Section \ref{discreteproc}\\
  $A^{n,k}$ & Time-change additive functional from  $X^n$ to $X^{n,k}$. & (\ref{ankdef})\\
  $\tau^{n,k}$ & Inverse of $A^{n,k}$. & (\ref{taudef}).\\
  $\ell^{n,k}$ & Occupation times for $J^{n,k}$. & (\ref{occtime})\\
  $L^{n,k}$ & Local times for $J^{n,k}$. & (\ref{loctime})\\
  $\hat{A}^{n,k}$ & Additive functional defined using $L^{n,k}$. & (\ref{tildeankdef})\\
  $\hat{\tau}^{n,k}$ & Inverse of $\hat{A}^{n,k}$. & Section \ref{discreteproc}\\
  $\hat{X}^{n,k}$ & Time-changed version of $J^{n,k}$ defined using $\hat{\tau}^{n,k}$. & (\ref{hatx}).\\
  $L^{(k)}$ & Local times of $B^{(k)}$. & Lemma \ref{localtimekcont}\\
  $\hat{A}$ & Additive functional defined using $L^{(k)}$. & (\ref{tildeakdef})\\
  $\breve{J}^{n,k},\breve{\ell}^{n,k},\breve{L}^{n,k}$ & Rescaled versions of   ${J}^{n,k},{\ell}^{n,k},{L}^{n,k}$ on $\breve{\underline{T}}_n(k)$. & Section \ref{finiteconverge}\\
  $B^{n,k}$ & Brownian motion on $(\breve{\underline{T}}_n(k), \lambda_n^{(k)})$. & Lemma \ref{jnkbk}\\
  $h^{n,k}$ & Hitting times of ``graph vertices'' by $B^{n,k}$. & (\ref{hitting})\\
  $\bar{L}^{n,k}$ & Local times of $B^{n,k}$. & Section \ref{finiteconverge}
\end{tabular}
\end{center}
Processes with a tilde represent the corresponding process mapped into $l^1$ by the distance-preserving maps $\psi$ or $\psi_n$ as appropriate. See Section \ref{crtsec} for the definition of $\psi$ and the end of Section \ref{cbptree} for the definition of $\psi_n$.

\section{Simple random walk estimates}

Collected in this section are several estimates for simple random walks on graphs, which are used in proving the convergence of the local times of jump processes on finite trees to those of the related Brownian motion, see Section \ref{finiteconverge}. We also prove results about the moments of the ``exit time'' of a simple random walk from a graph tree that are applied to prove the tightness result of Proposition \ref{anktildeank}.

\subsection{Occupation time tail bound for an interval}

We start by proving an exponential bound for the tail of the distribution of the occupation time of a simple random walk on an interval. In substance, the estimate was demonstrated by Kesten in the proof of \cite{Kestenunpub}, Lemma (4.64), but we include the proof here in order to state the result in a form that is more readily applicable in our situation.

We start by defining, for some fixed $R\in\mathbb{N}$, the sets $\Gamma_n:=\{0,1,\dots ,Rn\}$. Let $(Y^{n}_m)_{m\geq 0}$ be a discrete time simple random walk on $\Gamma_n$, starting from zero, under a probability measure $\mathbf{P}$. Here, we assume that vertices $x,y\in\Gamma_n$ are connected by an edge if and only if $|x-y|=1$. We also remark that the condition that $R$ is an integer is only for convenience, and the same argument can be applied for any $R>0$, when the interval considered is that between $0$ and $\lfloor Rn\rfloor$.

The processes $(\xi^n(x,m))_{m\geq 0,\:x\in\Gamma_n}$ will be the occupation time process for $Y^n$, i.e.
\[\xi^n(x,m):=\sum_{i=0}^{m}\mathbf{1}_{x}(Y^n_i).\]
The related hitting times will be written $(\varsigma^n_m(x))_{m\geq 1,\:x\in\Gamma_n}$, and can be defined by $\varsigma_m^n(x):=\min\{l:\:\xi^n(x,l)\geq m\}$. Finally, the intervals between the hitting times $(\sigma^n_m(x))_{m\geq 1,\:x\in\Gamma_n}$ are given by $\sigma_m^n(x):=\varsigma_{m+1}^n(x)-\varsigma_m^n(x)$. Note that, for fixed $x$ and $n$, $(\sigma^n_m(x))_{m\geq 1}$ is an independent, identically distributed sequence. We first prove a simple bound on the tail of the distribution of these intervals.

{\lem  Let $x\in\{0,1\}$. There exists $\varepsilon\in(0,1)$, $n_0\in\mathbb{N}$, depending only on $R$, such that
\begin{equation}\label{laha}
\mathbf{P}(\sigma_1^n(x)\geq \varepsilon n^2)>\frac{1}{4Rn},\hspace{20pt}\forall n\geq n_0.
\end{equation}}
\begin{proof} We will prove the result for $x=0$, the proof for $x=1$ is almost identical. By conditioning on the first step of the simple random walk, and then using the strong Markov property, we have
\begin{eqnarray*}
\lefteqn{\mathbf{P}(\sigma_1^n(x)\geq \varepsilon n^2)\geq \mathbf{P}(\mbox{$Y^n$ hits $Rn-1$ before returning to $0$}\:|\:Y_1^n=1)}\\
& & \hspace{80pt}\times \mathbf{P}(\mbox{$Y^n$ spends $\geq \varepsilon n^2$ in $[1,Rn-1]$ before hitting $0$}\:|\:Y^n_0=Rn-1).
\end{eqnarray*}
The first probability here is exactly equal to $(Rn-1)^{-1}$, by an elementary calculation. The final term is bounded below by
\[\mathbf{P}(\mbox{$S$ spends $\geq \varepsilon n^2$ in $[0,Rn-2]$ before hitting $Rn-1$}\:|\:S_0=0),\]
where $S=(S_n)_{n\geq 0}$ is a simple random walk on $\mathbb{Z}$. As $n\rightarrow\infty$, Donsker's theorem (see \cite{Kallenberg}, Theorem 14.9, for example) implies that the final term converges to
\[\mathbf{P}(\mbox{$B$ spends $\geq \varepsilon$ in $[0,R)$ before hitting $R$}\:|\:B_0=0),\]
where $B=(B_t)_{t\geq 0}$ is a standard Brownian motion on $\mathbb{R}$. Clearly, by taking $\varepsilon$ small, this probability can be chosen to be arbitrarily close to 1. The result follows.
\end{proof}

The independence of the sequences $(\sigma^n_m(x))_{m\geq 1}$ easily allows us to extend this result to the desired exponential bound.

{\lem \label{Occ} Let $x\in\{0,1\}$. There exist constants $c_{1}$ and $c_{2}$ depending only on $R$, such that
\[\mathbf{P}(\xi^n(x,n^2)\geq tn)\leq c_{1}e^{-c_{2}t},\hspace{20pt}\forall t\geq 0,\:n\in\mathbb{N}.\]}
\begin{proof} Let $x\in\{0,1\}$, $t\in\mathbb{N}$, and choose $\varepsilon$ and $n_0$ to satisfy the bound at (\ref{laha}). By definition, we have for $n\geq n_0$,
\begin{eqnarray*}
\mathbf{P}(\xi^n(x,n^2)\geq t n) &\leq & \mathbf{P}\left(\sum_{m=1}^{tn-1}\sigma^n_{m}(x)\leq n^2\right)\\
&\leq &\mathbf{P}\left(\sum_{m=1}^{t n-1}\mathbf{1}_{\{\sigma^n_{m}(x)\geq \varepsilon n^2\}}\leq \varepsilon^{-1}\right)\\
&\leq &\mathbf{P}\left(\mbox{Bin}(tn-1,\frac{1}{4Rn})\leq \varepsilon^{-1}\right)\\
&\leq & c_{3}e^{-c_{4}t},
\end{eqnarray*}
for some constants $c_{3}$ and $c_{4}$ depending only on $R$. Here, $\mbox{Bin}(n,p)$ represents a binomial random variable with parameters $n$ and $p$. Note also that we use the previous lemma for the third inequality. This estimate is easily extended to all $t$ and $n$ in the desired range by adjusting the constants suitably.
\end{proof}

\subsection{Crossing a tree}

Consider a graph tree $T$. Suppose that the shortest path between two vertices $x$ and $y$ in $T$ is of length $L$, and that the vertices $x$ and $y$ have degree $D_1$ and $D_2$ respectively. Assume that under the probability measure $\mathbf{P}$, the process $(X_m)_{m\geq 0}$ is a discrete time simple random walk on $T$ started from $x$. Denote by $N$ the number of visits by $X$ to $y$ before its first return to $x$. By observing the random walk on the path between $x$ and $y$, it is an elementary exercise to calculate that the exact distribution of $N$ is given by
\[\mathbf{P}(N=k)=\frac{1}{L^2D_1D_2}\left(1-\frac{1}{LD_2}\right)^{k-1},\]
for $k\geq 1$, and $\mathbf{P}(N=0)=1-1/LD_1$. These formulae allow us to deduce that if we define $\eta:=(N/D_2)-(1/D_1)$, then $\mathbb{E}\eta=0$; moreover, for every $k\geq 1$, there exists a constant $c=c(k)$, which does not depend on $D_1$, $D_2$ or $L$, such that
\begin{equation}\label{etamoments}
\mathbb{E}|\eta|^k\leq cL^{k-1}.
\end{equation}

\subsection{Exit times from a tree}

In the following lemma, $T$ is a rooted graph tree, with root $\rho$. The height of $T$ will be written $h(T)$. For a vertex $x\in T$, we write $T_x$ to denote the sub-tree determined by those vertices $y$ of $T$ such that the path from $\rho$ to $y$ passes through $x$. Suppose now that we join $D$ vertices to the root, each connected by a single edge. We shall denote by $\alpha(T,D)$ the expected time for a discrete time simple random walk on the graph consisting of $T$ and the $D$ extra vertices to leave the set of vertices of $T$, given that it started from $\rho$ (alternatively, this is the expected time for the walk to hit one of the extra vertices), and we shall denote by $\beta(T,D)$ the second moment of this time.

{\lem \label{alphabeta} For any graph tree $T$, and $D\geq 1$,
\[\alpha(T,D)=\frac{2|T|-2+D}{D},\hspace{20pt}\beta(T,D)\leq 36 (D+h(T))\frac{|T|^2}{D}.\]}
\begin{proof} The expression for $\alpha(T,D)$ is standard, see \cite{Kesten}, Lemma (2.28) for an example of a proof. In the same reference, it is also proved that
\[\beta(T,D)\leq 4\alpha(T,D)^2+ \frac{32h(T)}{D}\sum_{x\in T:\:x\sim \rho}|T_x|^2,\]
where $x\sim \rho$ means that $x$ is connected to $\rho$ by an edge. The sum is clearly bounded by $|T|^2$, and from the formula for $\alpha(T,d)$ we have that $\alpha(T,d)\leq 1+2|T|\leq 3|T|$. The result is easily deduced from these facts.
\end{proof}

\def\cprime{$'$}
\providecommand{\bysame}{\leavevmode\hbox to3em{\hrulefill}\thinspace}
\providecommand{\MR}{\relax\ifhmode\unskip\space\fi MR }
\providecommand{\MRhref}[2]{%
  \href{http://www.ams.org/mathscinet-getitem?mr=#1}{#2}
}
\providecommand{\href}[2]{#2}

\end{document}